\magnification=1200
\documentstyle{amsppt}

\global\newcount\numsec
\global\newcount\numfor
\global\newcount\numtheo
\global\advance\numtheo by 1

\def\senondefinito#1{\expandafter\ifx\csname#1\endcsname\relax}

\def\SIA #1,#2,#3 {\senondefinito{#1#2}%
\expandafter\xdef\csname #1#2\endcsname{#3}\else
\write16{???? ma #1,#2 e' gia' stato definito !!!!} \fi}

\def\etichetta(#1){(\veroparagrafo.\veraformula)%
\SIA e,#1,(\veroparagrafo.\veraformula) %
\global\advance\numfor by 1%
\write15{\string\FU (#1){\equ(#1)}}%
\write16{ EQ #1 ==> \equ(#1) }}

\def\letichetta(#1){\veroparagrafo.\verotheo
\SIA e,#1,{\veroparagrafo.\verotheo}
\global\advance\numtheo by 1
 \write15{\string\FU (#1){\equ(#1)}}
 \write16{ Sta \equ(#1) == #1 }}

\def\tetichetta(#1){\veroparagrafo.\veraformula 
\SIA e,#1,{(\veroparagrafo.\veraformula)}
\global\advance\numfor by 1
\write15{\string\FU (#1){\equ(#1)}}
\write16{ tag #1 ==> \equ(#1)}}

\def\FU(#1)#2{\SIA fu,#1,#2 }

\def\etichettaa(#1){(A\veroparagrafo.\veraformula)%
\SIA e,#1,(A\veroparagrafo.\veraformula) %
\global\advance\numfor by 1%
\write15{\string\FU (#1){\equ(#1)}}%
\write16{ EQ #1 ==> \equ(#1) }}

\def\BOZZA{
\def\alato(##1){%
 {\rlap{\kern-\hsize\kern-1.4truecm{$\scriptstyle##1$}}}}%
\def\aolado(##1){%
 {
{
 \rlap{\kern-1.4truecm{$\scriptstyle##1$}}}}}
}

\def\alato(#1){}
\def\aolado(#1){}

\def\veroparagrafo{\number\numsec}
\def\veraformula{\number\numfor}
\def\verotheo{\number\numtheo}

\def\Eq(#1){\eqno{\etichetta(#1)\alato(#1)}}
\def\eq(#1){\etichetta(#1)\alato(#1)}
\def\teq(#1){\tag{\aolado(#1)\tetichetta(#1)\alato(#1)}}
\def\Eqa(#1){\eqno{\etichettaa(#1)\alato(#1)}}
\def\eqa(#1){\etichettaa(#1)\alato(#1)}
\def\eqv(#1){\senondefinito{fu#1}$\clubsuit$#1
\write16{#1 non e' (ancora) definito}%
\else\csname fu#1\endcsname\fi}
\def\equ(#1){\senondefinito{e#1}\eqv(#1)\else\csname e#1\endcsname\fi}

\def\Lemma(#1){\aolado(#1)Lema \letichetta(#1)}%
\def\Theorem(#1){{\aolado(#1)Theorem \letichetta(#1)}}%
\def\Proposition(#1){\aolado(#1){Proposition \letichetta(#1)}}%
\def\Axioms(#1){\aolado(#1){Axioms \letichetta(#1)}}%
\def\Corollary(#1){{\aolado(#1)Corollary \letichetta(#1)}}%
\def\Remark(#1){\aolado(#1){Remark \letichetta(#1)}}%
\def\Definition(#1){{\noindent\aolado(#1){\bf Definition \letichetta(#1)$\!\!$\hskip-1.6truemm}}}
\def\Defprop(#1){\aolado(#1){Definition-Proposition \letichetta(#1)}}%
\def\Nada(#1){\aolado(#1){\letichetta(#1)}}%
\def\Example(#1){\aolado(#1) Example \letichetta(#1)$\!\!$\hskip-1.6truemm}

\def\include#1{
\openin13=#1.aux \ifeof13 \relax \else
\input #1.aux \closein13 \fi}

\openin14=\jobname.aux \ifeof14 \relax \else
\input \jobname.aux \closein14 \fi
\openout15=\jobname.aux

\input xypic

 
\pagewidth{16truecm}
\pageheight{22truecm}

\def\mapright#1{\smash{\mathop{\longrightarrow}\limits^{#1}}}
\def\dmapright#1#2{\smash{\mathop{\longrightarrow}\limits^{#1}_{#2}}}

\def\hr{\hookrightarrow}

\def\colim{\mathop{co\,lim}}

\def\Pi{P_i}

\def\q{\quad}
\def\C{\Cal C}

\def\E{\Cal E}
\def\A{\Cal A}
\def\Con{\Cal Con}

\def\Gf{\Gamma\!\!\!\Gamma_F}
\def\Lf{\Lambda\!\!\!\!\Lambda_F}
\def\GE{\Gamma\!\!\!\Gamma_E}
\def\Tf{\Theta_F}

\def\Ens{\Cal Ens}

\def\ProC{\Cal Pro(\Cal C)}

\def\EnsG{t {\Cal Ens}^{G}}
%
%

\topmatter
\title On the Galois Theory of Grothendieck\endtitle
\author by Eduardo J. Dubuc  and  Constanza S. de la Vega \endauthor
\endtopmatter

\vskip 4pt

\subheading{Spanish Abstract}

\vskip 4pt

 En este articulo tratamos la interpretacion hecha por Grothendieck
      de la teoria de Galois (y su relacion con el grupo fundamental y 
      la teoria de cubrimientos) en Expose V section 4, "Conditions 
      axiomatiques d'une theorie de Galois'' en el SGA1 1960/61.
      
           Esta es una hermosa muestra de matematicas muy rica en conceptos 
     categoricos, y su alcance es mucho mas vasto que el del trabajo
      original de Galois (asi como este llegaba mucho mas lejos que la 
      simple no resubilidad de la quintica).
      
           Aqui introducimos algunos axiomas y demostramos un teorema de 
     caracterizacion de la categoria (topos) de acciones de un grupo 
     discreto. Este teorema corresponde exactamente al resultado 
     fundamental de Galois. El teorema de Grothendieck caracteriza 
     la categoria (topos) de acciones continuas de un grupo topologico
      profinito. Desarrollamos una demostracion de este resultado 
      como un "paso al limite" (en un limite inverso de topos) de
       nuestro teorema de caracterizacion del topos de acciones de un
        grupo discreto. 
     Tratamos el limite inverso de topos trabajando con un simple 
     colimite filtrante (union) de las categorias que son sus 
     (respectivos) sitios de definicion. No se necesitan conocimientos
      avanzados de teoria de categorias para leer este articulo.

\vskip 4pt

\subheading{Introduction}

\vskip 4pt

In this paper we deal with Grothendieck's interpretation of Galois's 
Galois Theory (and its natural relation with the fundamental 
group and the theory of coverings) as he developed it in Expos\'e V,
section 4, ``Conditions axiomatiques d'une theorie de Galois'' 
in the SGA1 1960/61, [6]. 

This is a beautiful piece of mathematics very rich in categorical 
concepts, and goes much beyond the original Galois's scope (just as Galois 
went much further than the non resubility of the quintic equation). We show 
explicitly how Grothendieck's abstraction corresponds to Galois work.

 We introduce some axioms and prove a theorem of characterization
   of the category (topos) of actions of a discrete group. This theorem
   corresponds exactly to Galois fundamental result. The theorem of 
  Grothendieck characterizes the category (topos) of continuous actions
   of a profinite topological group. We develop a proof of this result
  as a "passage into the limit'' (in an inverse limit of topoi) of
  our theorem of characterization of the topos of actions of a 
  discrete group. 
  We deal with the inverse limit of topoi just working with an ordinary
  filtered colimit (or union) of the small categories which are their 
   (respective) sites of definition.

We do not consider generalizations of Grothendieck's work in [6], except by 
commenting briefly in the last section how to deal with the 
prodiscrete (not profinite) case. We also mention the work of Joyal-Tierney, 
which falls naturally in our discussion. 
   
    There is no need of advanced
   knowledge of category theory to read this paper, exept for the 
   comments in the last section..

\vskip 4pt

The paper has five sections, and a last one with comments and 
possible generalizations.

I. Examples.

II. Transitive actions of a discrete group.

III. Continuous transitive actions of a profinite group.

IV. All continuous actions of a profinite group.

V. All actions of a discrete group (and of a discrete monoid).

VI. Comments on this paper.

\vskip 6pt
     
\subheading{I - Examples}
\numsec=1\numtheo=1\numfor=1

\subheading{A - Classical Galois Theory}

\vskip 6pt

The ideas behind Galois Theory were developed through the work of Newton,
Lagrange, Galois, Kronecker, Artin and Grothendieck.

Before Galois it was known the following, that we write here using
modern notation:
%
%
   
\vskip 3pt

\proclaim{}
Given a field $k$ of numbers (of characteristic 0), a
polynomial of degree $n$, $f(x) \in k[x]$ with all roots
$\theta_1, \theta_2, \cdots, \theta_n$ different, and
a polynomial $s \in k[x_1,x_2,\cdots,x_n]$, then:

$$
s(\theta_1, \cdots, \theta_n)=s(\sigma\theta_1, \cdots, \sigma\theta_n)
\q \forall \sigma \in S_n \Rightarrow s(\theta_1, \cdots, \theta_n)
\in k
\Eq(1)
$$

where $S_n$ is the symmetric group in $n$ elements.

\endproclaim

\vskip 8pt

Galois in ``Memoire sur les conditions de resolubilit\'e des equations 
par radicaux'' (see [4]), proved a sharper result, he showed that:

\proclaim{}
There exists a group $G \subseteq S_n$ such that

$$
s(\theta_1, \cdots, \theta_n)=s(\sigma\theta_1, \cdots, \sigma\theta_n)
\q \forall \sigma \in G \Rightarrow s(\theta_1, \cdots, \theta_n)
\in k.
$$

\endproclaim

\vskip 6pt

He also showed the reverse implication.

Today this group is known as the Galois group $G(L/k)$ of
the splitting field $k \to L$ of $f$ over $k$, and the statement
above, takes the form:

Given any $\alpha \in L$, then:

$$
\sigma \alpha=\alpha \q \forall \sigma \in G(L/k) \Rightarrow
\alpha \in k.
$$

Notice that $\alpha=s(\theta_1, \cdots,\theta_n)$
for some $s \in k[x_1,x_2,\cdots,x_n]$.

See [4].

\vskip 4pt

The whole statement known as the Fundamental Theorem of Galois
Theory, after Artin, is the following:

\proclaim{}
Let $k \to N \to L$ be an intermediate extension of $L$
and $H \subseteq G(L/k)$ any subgroup. Then the assignments:

$$
\matrix
N&\longrightarrow&G(L/N)\subseteq G(L/k) \\
L^H&\longleftarrow&H\subseteq G(L/k)
\endmatrix
$$

where $G(L/N)=\{\sigma / \sigma \alpha=\alpha \q \forall \alpha \in N\}$
and $L^H=\{\alpha \in L / \sigma \alpha=\alpha \q \forall \sigma \in H\}$,
establish a one-one (contravariant) correspondence between the
lattices of subgroups of $G(L/k)$ and the subextensions of $L$.
\endproclaim

In particular, the Galois group $G(L/N)$ completely determines
the extension $N$.

\vskip 10pt

\subheading{B - Classical Theory of Coverings and the Fundamental Group}

\vskip 6pt

A local homeomorphism $X \mapright{p} B$ (sheaf over $B$) is a covering
space of $B$ if it is locally constant. That is, there exist a set
$S$, an (open) covering $U_\alpha \hookrightarrow B$ and homeomorphisms
$\varphi_\alpha$ such that:
%
%

$$
\diagram
\varphi_\alpha :U_\alpha \times S \drto_{\pi_{1}} \rto^>>>>>{\cong} 
& X_{|U_\alpha} \rto \dto &X \dto^{p} \\
&U_\alpha \rto|<\ahook &B 
\enddiagram
$$

\vskip 3pt

The set $S$ is the ``fiber'' of $p$. Given $b_0 \in U_{\alpha_0}$,
the homeomorphism $\varphi_{\alpha_0}$ establishes a bijection
$S \cong p^{-1}(b_0)$.

It is assumed that $B$ is path connected, locally path connected and
semi locally simply connected (see [11]) so that the existence of the
universal covering $\tilde{B} \to B$ is guaranteed.

The fundamental group $\pi_1(B,b_0)$ acts on the fiber
$S \cong p^{-1}(b_0)$.
We shall consider now only connected coverings which in this
context is equivalent to the fact that this action is transitive.

\vskip 4pt

The fundamental theorem of covering theory establishes:
\item{1)} The $\pi_1(B,b_0)$-set $S$ completely determines the
covering.
\item{2)} For any transitive $\pi_1(B,b_0)$-set $E$, there is
a covering with fiber $E$.

\vskip 3pt

The covering corresponding to $E$ is constructed as follows:

Take an $x \in E$ and let $H \subseteq \pi_1(B,b_0)$ be the
subgroup $H=\{g/ gx=x\}$. This subgroup $H$ acts on the fibers
of the universal covering and if we divide fiber by fiber by
this action we obtain a set $\tilde{B}/H$ that can be
given a topology making it a quotient space $\tilde{B} \to \tilde{B}/H$
which is a covering  $\tilde{B}/H \to B$.

\vskip 4pt

Recall that given any group $G$ there is a one-one correspondence
between the subgroups of $G$ (actually conjugacy classes of 
subgroups) and the transitive $G$-sets $E$ (actually equivalent 
classes of transitive $G$-sets under isomorphism).
Briefly, given $E$ take any $x \in E$ and set $H=Fix(x)$, on the other
hand, given $H$ set $E=G/H$.

Under this light, the fundamental theorem of coverings says the following:

\proclaim{}
Let 

$$
\diagram
\tilde{B} \rrto \drto && X \dlto \\
&B
\enddiagram
$$

be any covering of $B$ and $H\subseteq \pi_1(B,b_0)$ any subgroup.

\vskip 4pt

Let $Fix(x) \subseteq \pi_1(B,b_0)$, $Fix(x)=\{g/ gx=x \; \text{for the action
of} \; \pi_1(B,b_0) \; \text{on} \; S\}$ and 
$\tilde{B}/H$ be the quotient of $\tilde{B}$ by the action of $H$. 
Then the assignment:

$$
\matrix
X\to B&\longrightarrow&Fix(x) \q \text{where} \q x\in S \\
\\
\tilde{B}/H&\longleftarrow&H\subseteq \pi_1(B,b_0)
\endmatrix
$$

\vskip 4pt
%
%
establishes a one-one (covariant) correspondence between the
lattices of subgroups of $\pi_1(B,b_0)$ and the connected
coverings over $B$.

\endproclaim

\vskip 4pt

The above are two examples of a common theory that we call the
representable case of the Grothendieck axiomatic approach to
Galois Theory. This theory characterizes the category of all transitive
actions of a discrete (non necessarily finite) group.

\vskip 10pt

\subheading{II - Transitive actions of a discrete group}
\numsec=2\numtheo=1\numfor=1

\vskip 4pt

Let $\C$ be any category.

\vskip 4pt

\proclaim{\Definition (strictepi1)}

1) An arrow $X \mapright{f} Y$ in any category $\C$ is an strict 
epimorphism if given any compatible arrow $X \mapright{g} Z$, there 
exists a unique $Y \mapright{h} Z$ such that $g = h \circ f$. $g$ is 
compatible if for all 
$\C  \mapright{x} X$, $\C  \mapright{y} X$
with $g \circ x = g \circ y$, then also $f \circ x = f \circ y$.
\endproclaim

\vskip 4pt

It immediately follows that: $\text{strict epi  +  mono   =   iso}$

\vskip 4pt

Let now $A \in \C$ be any object, and let  
$G=(Aut(A))^{op}$ the opposite group of the group of automorphisms
of $A$.

\proclaim{\Definition (quotient)} Given any group $H$, 
a (left) action of $H$ on $A$ is a morphism of groups 
$H \to G=(Aut(A))^{op}$. The
categorical quotient $A \mapright{q} A/H$ is defined by the
following universal property:

\item{i)} $\forall h \in H, q\circ h=q$
\item{ii)} Given an arrow $A \mapright{x} X$ such that $\forall h \in H$,
$xh=x$, there exists a unique $A/H \mapright{\varphi} X$ such that
$\varphi \circ q=x$.

\vskip 3pt

In a diagram:

$$
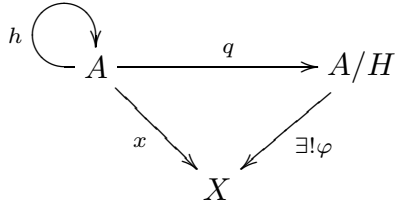
\diagram
\tolu^h A \drto_x \rrto^q && A/H \dlto^{\exists ! \varphi} \\
&X
\enddiagram
$$

We abuse the language and write `$h$' for the automorphism of $A$ 
corresponding to the element $h \in H$. Notice that different 
elements may define the same automorphism.

\endproclaim

Notice that it clearly follows that the arrow $A \mapright{q} A/H$ is 
an strict epimorphism. 

\proclaim{\Proposition (action)}
Given any $H\subseteq G$, $H$ acts on the left on the sets $[A,X]$
of morphisms in $\C$ from $A$ to any other object $X \in \C$:
\endproclaim

$$
\matrix
H&\times&[A,X]&\longrightarrow&[A,X] \\
h&&x&\longrightarrow&hx=x\circ h
\endmatrix
$$

where ``$\circ$'' is composition of arrows in $\C$.

\vskip 4pt

Assume for the moment that this action $G\times [A,X] \to [A,X]$
is transitive. We have in this way a functor:

\spreaddiagramcolumns {1pc}
$$ 
\diagram
\C \rto^{[A,-]_G} & \EnsG 
\enddiagram
$$
\spreaddiagramcolumns {-1pc}

where by $[A,X]_G$ we indicate the set $[A,X]$ with this
action of $G$ and $\EnsG$ indicates the category of transitive
$G$-sets.

\vskip 4pt

\vskip 4pt
%
%
\proclaim{\Remark (free)}
\endproclaim

The group $G$ has a canonical action on itself (given by left 
multiplication), thus we can consider
$G \in \EnsG$. $G$ furnished with this action is the free $G$-set
on one generator $e \in G$ ($e$ the neutral element). This means
that given any $E \in \EnsG$ there is a one-one correspondence:

$$
\aligned
& x\in E \over  \varphi:G\mapright{} E, \q \varphi(g)=gx \\
& \q \q e \mapsto x
\endaligned
$$

where the arrow $\varphi$ is a morphism of action (homogeneous map).

\vskip 6pt

Given any $x\in E$, the corresponding map $G \mapright{\varphi} E$
makes of $E$ a categorical quotient $E=G/H$ where $H=Fix(x)=
\{g\in G/ gx=x\}$.

\vskip 12pt

Consider any object $X \in \C$ and the $G$-set $[A,X]_G$. The
bijection in the remark above takes
the form:

$$
\aligned
{A\mapright{x} X \over  \varphi:G \mapright{} [A,X]_G} 
\endaligned
\Eq(bijection)
$$

This bijection means exactly that the object $A$ is the value
of the left adjoint of the functor $[A,-]_G$ evaluated in $G$. This
left adjoint is denoted:

\spreaddiagramcolumns {1pc}
$$
\diagram
\C & \EnsG \lto^{A\times_G (-)} 
\enddiagram
$$
\spreaddiagramcolumns {-1pc}

Thus $A\times_G G=A$.

\vskip 4pt

Consider now any $E \in t \Ens^{G}$ and an element $x_{0} \in E$.
Let $H=Fix(x_{0})$, and assume that the quotient $A/H$ exists in $\C$.
Then, it follows immediately that in the correspondence 
\equ(bijection) above, an
arrow $A \mapright{} X$ factors through $A/H$ if and only if the
corresponding arrow $G \to [A,X]_{G}$ factors through $G \to G/H \cong 
E$. Thus, there is a one to one correspondence:

$$
\aligned
{A/H \mapright{} X \over E \mapright{} [A,X]_G} 
\endaligned
$$

Therefore, the value of the left adjoint evaluated in the object $E$ is
equal to $A/H$:

$$ A\times_G E=A/H $$

Thus, if quotients of $A$ by subgroups of $Aut(A)$ exists in $\C$,
the functor $[A,-]_G$ has a left adjoint given by the formula above.

\vskip 8pt

We shall examine now conditions on $\C$ that will ensure that this
left adjoint together with $[A,-]_G$ establishes
an equivalence of categories.

\vskip 6pt
%
%
\proclaim{\Nada (axiomasRC)}{\bf Axioms for the representable connected case}
\endproclaim

\vskip 4pt

Consider a category $\C$ and an object $A \in \C$:

\vskip 4pt

\item{RC0)} For all X in $\C$ there exist $A \to X$, and every arrow $A 
\to X$ is a strict epimorphism (see \equ(strictepi1)).

\item{RC1)} 
The quotient $A \mapright{q} A/H$ of $A$ by any subgroup of $Aut(A)$,
$H\subseteq Aut(A)$, exists in $\C$ and it is preserved by the functor
$[A,-]$.

\item{RC2)} Every endomorphism of $A$ is an automorphism. That is
$[A,A]=Aut(A)$.

\vskip 4pt

Remark that the axiom RC2) follows (from \equ(PRC0)) if we assume that the 
functor $[A,-]$ preserves strict epimorphisms and that $[A,A]$ is a finite
set (see \equ(PC00) below).

\vskip 4pt

Notice that it follows immediately from RC0 that every arrow $X \to Y$ in $\C$ 
is an strict epimorphism. It also follows only from RC0 the following:

\vskip 3pt
\proclaim{\Proposition (PRC0)}
Axiom RC0) implies that the functor $[A, -]$ is faithful, reflects monomorphisms 
and reflects isomorphisms.
\endproclaim

\proclaim{Proof} \endproclaim
See \equ(PC0) where the statement is proved for any functor F in place 
of [A, -].

\vskip 4pt
\proclaim{\Remark (remarksRC)} \endproclaim 

{\bf .1 On RC1)}

The map $[A, A] \mapright{q_{*}} [A,A/H]$ factors $q_{*} = \eta \circ \rho$ 

$$
\diagram
[A,A] \rrtou^{q_{*}} \rto^{\rho} &[A,A]/H \rto^{\eta}&[A,A/H] 
\enddiagram
$$

RC1) means that $\eta$ is a bijection. This is divided in two parts:
     
\item{i)} if  $q \circ f = q \circ g$  then there exists  $h \in H$  such that 
 $f = h \circ  g$
 
\item{ii)} For all $A \mapright{x} A/H$, there exists $A \mapright{f} A$ 
such that $q \circ f = x$.

\vskip 3pt
{\bf .2 On RC2)}

In the presence of axiom RC2 condition i) above becomes equivalent to:

\item{iii)} $q \circ f = q$ implies $f \in H$

Also, given any arrow $A \mapright{x} X$, if we consider the surjective-injective 
factorization $x_{*} = \psi \circ \rho$:
$$
\diagram
[A,A] \rrtou^{x_{*}} \rto^{\rho} & I \rto^{\psi}&[A, X] 
\enddiagram
$$

axiom RC2 implies that $I = [A, A]/H$, where $H = Fix(x) \subset Aut(A)$
  
\vskip 4pt
%
%
In the next three propositions of this section all the three axioms are 
needed.

\proclaim{\Proposition (PRC3)}
Any arrow $A \mapright{x} X$ is a quotient of $A$ by the 
subgroup $H \subseteq Aut(A)$, $H=Fix(x)$ \endproclaim

\proclaim{Proof} \endproclaim

Consider the factorization:
   
$$
\diagram
A \rrtou^{x} \rto^{q} &A/H \rto^{\varepsilon} &X 
\enddiagram
$$

If we apply the functor $[A,-]$ to this diagram, we have:

$$
\diagram
[A,A] \drto_{\rho} \rrtou^{x_{*}} \rto^{q_{*}} &[A,A/H] \rto^{\varepsilon_{*}}&[A,X] \\
&[A,A]/H \uto^{\eta} \urto^{\psi} 
\enddiagram
$$

\vskip 3pt

Where $\rho$, $\eta$ and $\psi$ are as in \equ(remarksRC).1 and 
\equ(remarksRC).2.

Since $\eta$ is a bijection and $\psi$ is injective, it follows that  
$\varepsilon_{*}$ is also injective. 
Thus by proposition \equ(PRC0), $\varepsilon$ is a monomorphism. Since 
every arrow is an strict epimorphism, it follows that   
it is an isomorphism (see \equ(strictepi)).

\vskip 4pt

The statement in this proposition is actually equivalent to the fact that 
the functor $[A,-]$ reflects isomorphisms. In fact, given any arrow 
$X \mapright{f} Y$, take $A \mapright{x} X$. 
We have $Fix(x) \subset Fix(f \circ x)$. If $f_{*}$ is a bijection, given any 
$A \mapright{h} A$, the equation $f \circ x \circ h = f\circ x$ implies 
$x \circ h = x$. Thus $Fix(x) = Fix(f \circ x)$. 
It follows from the statement of proposition \equ(PRC3) that $f$ is an 
isomorphism.

\proclaim{\Proposition (pepis)} 
The functor $[A,-]$ preserves strict epimorphisms.
\endproclaim

\proclaim{Proof} \endproclaim
By RC0) it clearly suffices to consider the case $A \mapright{x} X$. 
Then the statement follows immediately from proposition \equ(PRC3).

\proclaim{\Proposition (PTA)}
The action of $Aut(A)$ on $[A,X]$ is transitive for all $X$.
\endproclaim

\proclaim{Proof} \endproclaim
Since every arrow $A \to X$ is an strict epimorphism (RC0), proposition 
\equ(pepis) means exactly that the action of the monoid $[A,A]$
is transitive. Then RC2) finishes the proof. 

\proclaim{\Theorem (adjoint1)} Let $\C$ be any category
and $A \in \C$ as above. If the axioms RC0) to RC2) hold then the 
left adjoint $A\times_G (-) \dashv [A,-]_G$ exists and the maps:

\vskip 2pt

\item{a)} $E\cong [A,A]/H \mapright{\eta} [A,A/H]$
%
%

\item{b)} $A/H \mapright{\varepsilon} X$

are isomorphisms. Thus, they establish an equivalence of categories:

\spreaddiagramcolumns {1pc}
$$
\diagram
\C  \rto<1ex>^{[A,-]_G} & t \Ens^G \lto<1ex>^{A \times_G (-)}  \\
\enddiagram
$$
\spreaddiagramcolumns {-1pc}

\endproclaim

Notice that $A/H=A\times_G [A,X]_G$ and $[A,A/H]=[A,A\times_G E]_G$,
and that $\eta$ and $\varepsilon$ are the unit and counit of the
adjunction.

\vskip 3pt

\proclaim{Proof} \endproclaim
That $\eta$ and $\varepsilon$ are isomorphisms is proved in 
\equ(remarksRC).1 and \equ(PRC3) respectively.

\vskip 10pt

It is not difficult but \underbar {not trivial} to check that the classical
case of Galois Theory (A in section I) and the covering example
(B in section I) satisfy axioms RC0) to RC2). 
In the first case $\C$ is the category dual to the category of subextension of 
the splitting field $L$ of $f$ and the object $A$ is the extension 
$L$, and in the second case $\C$ is the category of path connected
coverings of the space $B$ and the object $A$ is the universal covering.
The reader can also check that the recipes a),b) in theorem 
\equ(adjoint1) are 
also exactly the ones given in the examples. Notice that in example A, the 
subextension $L^H=\{ \alpha \in L/ \sigma \alpha=\alpha \q \forall 
\alpha \in H\}$ is precisely the quotient of $L$ by $H$ in the dual
category. 
 
\vskip 8pt

There are two conspicuous non-examples of the above theorem:

First take $k \to \tilde{k}$ be the algebraic closure of $k$
(we assume for simplicity that $k$ is of characteristic 0),
and $\C$ the dual category of all intermediate extensions.
In this case it is known that the theorem is false, thus some
of the axioms must not hold.

The other non example is the case like above of all path connected
coverings of $B$, a topological space path connected, locally path
connected but not necessarily semilocally simply connected. In this
case there is no universal covering and we don't have the object
$A \in \C$.

These two non-examples can be suitable modified so that a similar
theorem will hold. In the first case there have to be taken only
finite extension, and in the second case only coverings with finite
fibers. So, the algebraic closure or the universal covering, even 
when the latter exists, are not any more objects in the category.

Both these cases share the fact that there is no universal object $A$ 
in the sense of axiom RC0 in the category $\C$. 
Grothendieck in [6] replaces the representable functor
$[A,-]$ by a functor $F:\C \to \Ens$, which is a filtered colimit
of representables and establishes a theorem which has as particular
cases the two situations described above. 
In addition, he assumes that the functor $F$ takes its values on 
\underbar{finite} sets. Due to this assumption, his theorem is not a 
generalization of the theorem above, and in particular it  
\underbar{does not} has as a particular case the example B. Now, the group is
not longer the usual fundamental group of loops, but it is its profinite
closure, and the coverings included in the theorem 
%
%
are only those with finite fibers.
Nevertheless, it can also be proved a theorem that includes the case 
of all covering spaces, even when there is no universal covering, 
as it is said (but not proved) in [[1], 2.7.4 and 2.7.5] (see 
\equ(IIIyIV)).

\vskip 10pt

\subheading{III - Continuous transitive actions of a profinite group}
\numsec=3\numtheo=1\numfor=1

\vskip 6pt

Consider a category $\C$ and a functor $F:\C \to \Ens$.

Recall that the diagram of $F$, that we denote $\Gf$
is the category whose objects are
pairs $(a,A)$ where $a \in F(A)$ and whose arrows
$(a,A) \mapright{f} (a',A')$ are maps $A\mapright{f} A'$ such that
$F(f)(a)=a'$.

There is a functor

\spreaddiagramrows {-2pc}
$$
\diagram
\Gf \rto & \Ens^C \\
(a,A) \rto & [A,-] 
\enddiagram
$$
\spreaddiagramrows {2pc}

with the obvious definition on arrows, and $F$ is the colimit of
this diagram, that is:

$$F=\dsize\lim_{(a,A) \in \Gf}{[A,-]}$$

Recall also that there is a one-one correspondence:

$$
\aligned
{a \in F(A)  \over [A,-] \dmapright{}{a} F} 
\endaligned
$$

\vskip 4pt

where the arrow $a$ is the natural transformation completely 
characterized by the equation $a(id_A)=a$ (notice the abuse of 
notation).

\vskip 10pt

\proclaim{\Nada (axiomasC)} {\bf Axioms for the connected case}
\endproclaim

Consider the following axioms:

\item{C0)} For all $A \in \C, F(A) \neq \emptyset$  and
every arrow $A \mapright{x} B$ is an strict epimorphism (in this case 
the category $\Gf$ is actually a poset, see remark \equ(poset)).

\item{C1)} Given any object $A$ and a finite group
$H$ acting on $A$, $H \to (Aut(A))^{op}$, the quotient $A \mapright{q} A/H$
exists and it is preserved by $F$ (see definition \equ(quotient)).

\item{C2)} $F(A)$ is finite $\forall A \in \C$ and F preserves strict 
epimorphisms.

\item{C3)} The poset $\Gf$ has all finite meets, in particular it is 
filtered.

\proclaim{\Remark (poset)}
Since all maps $A \mapright{x} B$ in $\C$ are epimorphisms,
for any object $X \in \C$ the transition morphisms corresponding
to an arrow $(a,A) \mapright{x} (b,B)$ in $\Gf$,
$x^*:[B,X] \mapright{} [A,X]$ are all injective functions.
By construction of filtered colimits in $\Ens$ it follows
that the canonical maps of the colimit $[A,-] \mapright{a} F$
are injective natural transformations (thus monomorphisms
in the category $\Ens^\C$). This implies that $\Gf$ is a poset.
\endproclaim
%
%

\proclaim{\Remark (PC00)}
Axiom C2 implies that for all $A \in \C$ every endomorphism of $A$ is an 
automorphism. That is $[A,A]=Aut(A)$.
\endproclaim

\proclaim{Proof} \endproclaim
Given $A \mapright{h} A$, the function $F(A) \mapright{h_*} F(A)$ 
is surjective, then if $F(A)$ is finite it is a bijection, and thus
by \equ(PC0) $h$ is an isomorphism.

\proclaim{\Remark (OnC2)} 
Axiom C1) means that the factorization of the map $F(q)$ through $F(A)/H$
is a bijection. In a diagram:

\spreaddiagramcolumns {-0,7pc}
$$
\diagram
F(A) \drto \rrto^{F(q)} && F(A/H) \dlto^{\cong}  \\
&F(A)/H
\enddiagram
$$
\spreaddiagramcolumns {0,7pc}

\endproclaim

\vskip 4pt

\proclaim{\Remark (meets)}
Notice that the condition of having all finite meets on the poset 
$\Gf$ clearly implies it is filtered.
Filterness means for posets that given any pair $[A,-] \mapright{a} 
F$, $[B,-] \mapright{b} F$, there exist $[C,-] \mapright{c} F$, 
$C \mapright{x} A$ and $C \mapright{y} B$ such that

$$
\diagram
A &&[A,-] \drto_{x^*} \drrto^a \\
&C \dlto_y \ulto^x &&[C,-] \rto^c &F \\
B &&[B,-] \urto^{y^*} \urrto_b
\enddiagram
$$

commutes.

Existence of meets means that there is a last such $C$. In this
context, this is a strictly stronger condition, which is actually
needed in proposition \equ(closure).
\endproclaim

\vskip 4pt

Remark \equ(meets) also says
that $F$ is an strictly pro-representable functor in
the sense of Grothendieck [1], and its formal dual is a pro-object
that we call $P$. We have the following diagram:

$$
\diagram
\C \q \drto|<<\ahook \rrto|<\ahook && (\Ens^\C)^{op}  \\
&\ProC \urto 
\enddiagram
$$
%
%

where all the arrows are full and faithful functors, and to
$F \in (\Ens^\C)$  it correspond $P \in \ProC$.
To the colimit diagram $[A,-] \mapright{a} F$ (indexed by
$\Gf$) in $(\Ens^\C)$ it correspond a limit diagram
(indexed by $\Gf$) $P \mapright{a} A$ in $\ProC$.
Here by abuse, given an object $A \in \C$ we also denote by $A$ the
corresponding (representable) pro-object  in $\ProC$, $A \in \ProC$.
Thus, by this trick of Grothendieck, the functor $F$ becomes
representable, but from the outside of $\C$, by
a pro-object $P$. Notice that tautologically we have $F(X)=[P,X]$.

Notice also that (tautologically) the diagram of $F$, $\Gf$ becomes
the category of all objects of $\C$ below $P$: $P \mapright{a} A$,
with morphisms $A \mapright{x} B$ in $\C$ such that the triangle

$$
\diagram
P \dto_a \rto^b &B \\
A \urto_x 
\enddiagram
$$

commutes.
Furthermore given any object $A \in \C$ all maps
$P \mapright{a} A$ are epimorphisms in the category $\ProC$
(since they are monomorphisms in $\Ens^{\C}$).

\vskip 4pt

\proclaim{\Proposition (PC0)}
From axiom C0 it follows that the functor $F$ is faithful, reflects 
monomorphisms and reflects isomorphisms.
\endproclaim

\proclaim{Proof} \endproclaim
Let $A \mapright{f} B$, $A \mapright{g} B$ be such that $F(f)=F(g)$.
Take $a \in F(A)$ and let $b=F(f)(a)=F(g)(a)$. We have a commutative
diagram:

$$
\diagram
&P \dlto_{a} \drto^{b} \\
A \rrto<1ex>^{f} \rrto<-1ex>_{g}&&B
\enddiagram
$$

Since $a$ is an epimorphism, it follows $f=g$.
This shows the first assertion.

Now, let $A \mapright{f} B$ be such that $F(f)$ is an injective 
function,
and let $C \mapright{x} A$, $C \mapright{y} A$ be any two morphisms
such that $f \circ x=f \circ y$. Clearly  $F(f) \circ F(x) = F(f) \circ 
F(y)$. Thus $F(x)=F(y)$, and by the first assertion $x=y$.
This shows the second assertion. Clearly C0 implies that the third 
assertion follows from the second (see \equ(strictepi).

\vskip 4pt

\proclaim{\Definition (Galois)} An object below
$P \mapright{a} A$ is Galois $\Leftrightarrow$ the arrow
$Aut(A) \mapright{a^*} [P,A]$ is a bijection.
That is:

$$
\diagram
P \dto_a \rto^{\forall y} &A \\
A \urto_{\exists ! \approx h} 
\enddiagram
$$

\endproclaim
%
%

\vskip 4pt

From \equ(PC00) this is equivalent to the fact that the arrow
$[A,A] \mapright{a^*} [P,A]$ is a bijection.

\vskip 4pt
 
Notice that it follows immediately that $P \mapright{a} A$ 
is Galois $\Leftrightarrow$ $P\mapright{b} A$ is Galois $\forall b$. We
are justified in saying then that the object $A$ is Galois.

\vskip 4pt

Referring to section I, Galois objects are exactly the Galois extensions 
in example A and the regular coverings in example B.

\vskip 4pt

\proclaim{\Definition (CA)} Given any $P \mapright{a} A$
Galois, we define the full subcategory $\C_A \hookrightarrow \C$
as follows:

$X \in \C_A \Leftrightarrow [A,X] \mapright{a^*} [P,X]=F(X)$
is a bijection, that is:

$$
\diagram
P \drto_{\forall x} \rrto^{a} && A \dlto^{\exists !}  \\
&X 
\enddiagram
$$
\endproclaim

Observe that $A\in \C_A$ since $A$ is Galois and that $\C_A$ does not
depend on $a$ but only on the object $A$.

\proclaim{\Proposition (CB cont CA)} Given any transition morphism 
between Galois objects in $\Gf$,

$$
\diagram
&P \dlto_{a} \drto^{b}  \\
A \rrto^x &&B 
\enddiagram
$$

there is an inclusion of full subcategories $\C_B \hookrightarrow \C_A$.
\endproclaim

\vskip 6pt

\proclaim{\Remark (F=A)} The restriction of the functor
$F$ on $\C_A$ is the representable functor by $A$, that is

$$
\diagram
\C_A \drto_{[A,-]} \rto|<\ahook & \C \dto^{[P,-]=F}\\
&\Ens
\enddiagram
$$

commutes.
\endproclaim

\proclaim{\Proposition (CA=below A)} 
The category $\C_A$ is closed by quotients of actions of finite groups.
\endproclaim

\proclaim{Proof} \endproclaim

Suppose we have a quotient $A \mapright{q} A/H$.
Consider the commutative diagram:
%
%

$$
\diagram
[A,A] \dto_{q_*} \rto^{a^*} & [P,A] \dto^{q_*}  \\
[A,A/H] \rto_{a^*} & [P,A/H]
\enddiagram
$$

\vskip 6pt

We have to show that the lower $a^*$ arrow is a bijection. 
This is immediate since the upper $a^*$ is
bijective and the right $q_*$ is surjective since
$[P,A/H] \cong [P,A]/H$ by axiom C1, see \equ(OnC2).

\proclaim{\Theorem (equivalence)}
The category $\C_A$ with the object $A \in \C_A$ satisfies the axioms
RC0) to RC2). Thus by theorem \equ(adjoint1) we have an equivalence 
of categories:

\spreaddiagramcolumns {1pc}
$$
\diagram
\C_A  \rto<1ex>^{[A,-]_G} & t \Ens^G \lto<1ex>^{A \times_G (-)}  \\
\enddiagram
$$
\spreaddiagramcolumns {-1pc}

where $G=Aut(A)^{op}$.

\endproclaim

\proclaim{Proof} \endproclaim
Consider the remark \equ(F=A), then:
axiom C0) gives axiom RC0), axiom C1) together with proposition \equ(CA=below A) 
give axiom RC1) and proposition \equ(PC00) give axiom RC2).

\vskip 10pt

\proclaim{\Proposition (closure)} (Existence of Galois
closure) Given any $X \in \C$, there exist a Galois object
$P \mapright{a} A$ such that $\forall P \mapright{x} X$,
there exist a unique $\pi_x$ that makes the following diagram
commutative:

$$
\diagram
P \drto_{\forall x} \rrto^{a} && A \dlto^{\exists ! \pi_x}  \\
&X 
\enddiagram
$$

In other words, the map $a^*:[A,X] \to [P,X]$ is bijective.

\endproclaim

\proclaim{Proof} \endproclaim

Consider the finite set of all $P \mapright{x} X$ and let 

$$
\diagram
P \drto_{x} \rrto^{a} && A \dlto^{\pi_x}  \\
&X 
\enddiagram
$$ 
%
%

be the infima in $\Gf$ which exists by C3).
We shall see that $A$ is Galois.

Given any $P \mapright{y} A$, consider the map $[P,X] \to [P,X]$,
that sends $x \to \pi_x \circ y$.
This map is an injective map ($\pi_x \circ y=\pi_{x'} \circ y
\Rightarrow \pi_x=\pi_{x'} \Rightarrow x=x'$) between finite sets.
Thus it is a bijection. Given any $x$, $x$ is then of the form
$x=\pi_{x'} \circ y$ for some $x'$.
Consider the diagram:

$$
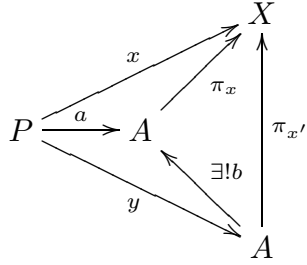
\diagram
&&X \\
P \urrto^x \rto^{a} \drrto_y &A \urto_{\pi_x}  \\
&& A \ulto_{\exists ! b} \uuto_{\pi_{x'}} 
\enddiagram
$$

\vskip 6pt

The unique existence indicated in the diagram follows by the
universal property of the infima. This shows that $A$ is
Galois.

\vskip 4pt

\proclaim{\Theorem (cofinal)} (Corollary of \equ(closure)) 
Galois objects are cofinal in the diagram of $F$. 
Thus, $F$ is a filtered colimit of representables $[A,-]$ with
$A$ Galois. Let $\Lf$ be the cofinal subdiagram of $\Gf$ whose
objects are Galois. We have $P \mapright{a} A$, $(a,A) \in \Lf$
is a filtered inverse limit diagram in such a way that $\C$ becomes a 
filtered colimit of the full subcategories $\C_{A} \hookrightarrow \C$.
\endproclaim

\proclaim{Proof} \endproclaim
It follows straightforward from proposition \equ(closure).

\vskip 15pt

Consider $(a,A) \mapright{x} (b,B)$ in $\Lf$. We define
$\rho_x:[A,A] \to [B,B]$ by means of the following commutative
diagram:

$$
\diagram
[P,A] \rto^{x_*} & [P,B] \\
[A,A] \uto_{a^*}^\cong \rto_{\rho_x} & [B,B] \uto_{b^*}^\cong 
\enddiagram
$$

\vskip 6pt

We define also $[P,P] \mapright{\pi_a} [A,A]$ by the following
commutative diagram:

$$
\diagram
[P,P] \drto_{a_*} \rrto^{\pi_a} && [A,A] \dlto^{a^*}  \\
&[P,A] 
\enddiagram
$$

\vskip 5pt

It follows that we have a cone over $\Lf$:

$$
\diagram
[P,P] \drto_{\pi_b} \rrto^{\pi_a} && [A,A] \dlto^{\rho_x}  \\
&[B,B] 
\enddiagram
$$

\vskip 6pt

and the diagram $[P,P] \mapright{\pi_a} [A,A]$ is an inverse
limit diagram of finite groups with the transition morphisms
$[A,A] \mapright{\rho_x} [B,B]$ surjective.
It follows immediately that $Aut(P)=[P,P]$.
The fact that $[P,P]$ is an inverse limit of finite groups means
by definition that it is a profinite group.
Moreover, the projections in the limit diagram
$[P,P] \mapright{\pi_a} [A,A]$ are also surjective (see 1) in 
proposition \equ(Glambda)).

\proclaim{\Definition (profinite group)} Let  $\Gamma \to Gr$
be an inverse filtered limit of finite groups.
The inverse limit group $\pi \mapright{\rho_\lambda} G_\lambda$ is
called a profinite group when considered as a topological group
with the product topology of the discrete finite groups $G_\lambda$.
\endproclaim

If $\pi$ acts continuously and transitively on a set (discrete space) 
$E$, given any $x \in E$ there is a continuous surjective function 
$\pi \mapright{(-) \cdot x} E$. Thus, since $\pi$ is compact, $E$ must 
be finite. But more than that is true:

\proclaim{\Proposition (cont=fact)} A transitive action 
$\pi \times E \to E$ on a set (discrete space) $E$ is continuous
if and only if it factors:

$$
\diagram
\pi \times E \dto_{\rho_{\lambda} \times id} \rto &E \\
G_\lambda \times E \urto 
\enddiagram
$$

for some $\lambda$. In particular $E$ must be a finite set.

\endproclaim

\proclaim{Proof} \endproclaim
In fact, consider any element $x \in E$, and let $H \subset \pi$ be 
the subgroup $H = Fix(x)$. $H$ is open because $E$ is discrete. Since 
the subgroups $K_{\lambda}=Ker(\rho_{\lambda})$ are a neighborhood base 
of the unit $e$, there is a $K_{\lambda} \subset H$. Since 
$K_{\lambda}$ is normal, it follows that $K_{\lambda} \subset Fix(y)$ for all 
other elements $y \in E$. It follows that the given action 
factors through $G_{\lambda}$.


\vskip 4pt

\proclaim{\Proposition (Glambda)} Let 
$\pi \mapright{\rho_\lambda} G_\lambda$ be a profinite group. Recall
that the
transition morphisms $G_\lambda \mapright{\rho_{\lambda \mu}} G_\mu$
(for any $\lambda \to \mu$ in $\Gamma$) are surjective.  We have:

\item{1)} The projections $\pi \mapright{\rho_\lambda} G_\lambda$
are surjective.

\item{2)} Given any $\lambda \in \Gamma$, there is a full and
faithful inclusion of categories:

$$
\diagram
t\Ens^{G_\lambda} \rto|<\ahook & ct\Ens^{\pi} \\
\enddiagram
$$

where the ``$t$'' stands for transitive and the ``$c$'' stands for
continuous.

\item{3)} Given any arrow $\lambda \to \mu$ in $\Gamma$, there is
a full and faithful inclusion of categories:

$$
\diagram
t\Ens^{G_\mu} \rto|<\ahook & t\Ens^{G_\lambda} \\
\enddiagram
$$

\item{4)} The category $ct \Ens^\pi$ is a filtered colimit (indexed 
by $\Gamma$) of the subcategories $t \Ens^{G_\lambda}$.

\endproclaim

\proclaim{Proof} \endproclaim

The statement 1) is by no means easy to proof, but it is 
systematically assumed to be true in the literature. We have not 
found yet a proof in print to give as a reference. The rest of the 
statements 2), 3) and 4) are straightforward.

\vskip 4pt
We establish now the compatibilities arising from a transition
morphism in $\Gf$:

$$
\diagram
&P \drto^{b} \dlto_{a} \\
A \rrto_{x} &&B 
\enddiagram
$$

\proclaim{\Proposition (commutes1)} The following diagram
commutes (up to a natural isomorphism):

\spreaddiagramcolumns {1pc}
$$
\diagram
\C_A \rto^{[A,-]_G} & t \Ens^G \\
\C_B \rto^{[B,-]_L} \uto|<\ahook  & t \Ens^L \uto|<\ahook 
\enddiagram
$$
\spreaddiagramcolumns {-1pc}

\vskip 6pt

where $G=Aut(A)^{op}$ and $L=Aut(B)^{op}$.
\endproclaim

\proclaim{Proof} \endproclaim

Let $X\in \C_B$ and consider the morphism in $t \Ens^G$,
$[B,X] \mapright{x^*} [A,X]$ induced by $x$. Is straightforward
to check that it is actually a morphism of $G$-actions (where
$[B,X]$ is considered with the action induced by $G \mapright{\rho_x} L$).
On the other hand, the diagram:

$$
\diagram
[A,X]  \rto^{a^*} &[P.X] \\
[B,X] \uto^{x^*} \urto_{b^*} 
\enddiagram
$$

clearly commutes.

Since $a^*$ and $b^*$ are bijections, $x^*$ is a bijection.

\proclaim{\Proposition (commutes2)}
The following diagram commutes (up to natural isomorphism):

\spreaddiagramcolumns {1pc}
$$
\diagram
\C_A & t \Ens^G \lto^{A\times_G (-)} \\
\C_B \uto|<\ahook  & t \Ens^L \lto^{B\times_L (-)} \uto|<\ahook 
\enddiagram
$$
\spreaddiagramcolumns {-1pc}

\endproclaim

\proclaim{Proof} \endproclaim

Since every object in $t \Ens^L$ is a quotient of $L$ (\equ(free)), and the
functors in the diagram preserves quotients 
(the horizontal ones are left adjoints, for the vertical inclusions 
use \equ(CA=below A)) it will be enough to show the proposition
for the object $L$.
Recall that given $A \mapright{x} B$,
then $B\cong A/H$ where $H \hookrightarrow G$, $H=Fix(x)$.
It is immediate to check that $G \to L$ is the quotient $L \cong G/H$.
With this in hand, it is clear that the object $L$ in $t \Ens^L$
goes into $B$ in $\C_A$ by any of the two paths in the diagram.

\vskip 4pt

\proclaim{\Proposition (action2)} Given any $X \in \C$,
there is a continuous transitive left action of the group 
$\pi=[P,P]^{op}$ on the set $[P,X]$:

$$
\matrix
\pi&\times&[P,X]&\mapright{}&[P,X] \\
h&\times&x&\mapright{}&x \circ h
\endmatrix
$$

given by composition in $\ProC$.

\endproclaim

\proclaim{Proof} \endproclaim

Let $P \mapright{a} A$ be such that $a^*:[A,X] \to [P,X]$ is a
bijection (given by proposition \equ(closure)).
Consider the following commutative diagram:

$$
\diagram
P \rto^{h} \dto_{a} & P \rto^{x} \dto_{a} & X\\
A \rto_{\tilde h} & A \urto_{\tilde x} 
\enddiagram
$$

where $\tilde{x}=(a^*)^{-1}(x)$ and $\tilde{h}=\pi_a(h)$
(see definition of $\pi_a$).
Thus, we have an action $x\circ h=\tilde{x}\circ \tilde{h} \circ a$.

Let $G=[A,A]^{op}$ and consider the following diagram:

$$
\diagram
\pi \times [P,X] \rto \dto_{\pi_a \times id} & [P,X] \ddto^{(a^*)^{-1}}\\
G \times [P,X] \dto_{id \times (a^*)^{-1}} \\
G \times [A,X] \rto & [A,X]  
\enddiagram
$$

\vskip 6pt

It only remains to see that this action is transitive.
That the action $G \times [A,X] \rightarrow [A,X]$ is transitive has 
been proved in proposition \equ(PTA). Now the statement follows
immediate since the morphism $\pi_a:\pi \to G$ is surjective.

\vskip 10pt

For any $X \in \C$, let $[P,X]_{\pi}$ be the
set $[P,X]$ enriched with this action.
In this way we have a functor:

\spreaddiagramcolumns {1pc}
$$
\diagram
\C \rto^{[P,-]_\pi} & ct\Ens^{\pi} \\
\enddiagram
$$
\spreaddiagramcolumns {-1pc}

into the category of continuous transitive $\pi$-actions.

We shall show now how from our theorem for the
representable case and the structure results 
in proposition \equ(CB cont CA) and proposition \equ(Glambda), 
it follows that the functor $[P,-]_{\pi}$ has
a left adjoint and establishes an equivalence of categories.

\proclaim{\Theorem (adjoint2)}
The functor 

\spreaddiagramcolumns {1pc}
$$
\diagram
\C \rto^{[P,-]_\pi} & ct\Ens^{\pi}
\enddiagram
$$ 
\spreaddiagramcolumns {-1pc}

has a left adjoint
$P\times_\pi (-)$ which establishes an equivalence of categories.
\endproclaim

\proclaim{Proof} \endproclaim

We have filtered colimits of full subcategories:

\spreaddiagramcolumns {1pc}
$$
\diagram
\C \rto^{[P,-]_\pi}& ct \Ens^\pi  \\
\C_A \uto|<\ahook \rto<1ex>^{[A,-]_G} & t \Ens^G \lto<1ex>^{A \times_G (-)} 
\uto|<\ahook 
\enddiagram
$$
\spreaddiagramcolumns {-1pc}

\vskip 6pt

\vfill \eject

which are compatible in the sense of propositions \equ(commutes1) and 
\equ(commutes2).
Given $E \in  ct \Ens^{\pi}$, to define $P\times_\pi E$, observe
that $E \in t \Ens^G $ for some $G=Aut(A)^{op}$ and define
$P\times_\pi E=A\times_G E$.
That this is well defined and actually determines a left adjoint
of $[P,-]_\pi$ follows immediately from the fact that the subcategories
are full and the quoted compatibilities.

Finally, since in all $A$-levels we have an equivalence, it follows
that we also have an equivalence in the limit.

\vskip 10pt

\subheading{IV - All continuous actions of a profinite group}
\numsec=4\numtheo=1\numfor=1

\vskip 6pt

We consider now Grothendieck axioms as he wrote them in [6] and prove 
its fundamental theorem. He considered a category with, in particular, 
finite sums, and proved that it is equivalent to the category of 
finite continuous actions on a profinite group.
However, in our proof of this result it can be seen clearly that the 
argument goes through if one assume arbitrary sums, and the conclusion 
is now that the category is equivalent to the category (in this case a 
topos) of all continuous actions of the same profinite group.
We shall write this two results in parallel: 

\proclaim{Convention}
The word [finite] between brackets will mean that the 
statement stands in fact for two statements: one assuming finite 
and the other not assuming finite.
When the word finite appears not in between brackets it has its usual 
meaning, and there is only one statement (which assumes finite).
\endproclaim

\vskip 6pt

\proclaim{\Nada (axiomasG)} {\bf Grothendieck's axioms}
\endproclaim

\vskip 4pt

Consider a category $\C$ and a functor $F:\C \to \Ens$ and assume the
following axioms:

\vskip 4pt

First we introduce the axiom G0) necessary to deal with the not 
finite case:

An object $X \in \C$ is called finite if $F(X)$ is a finite set.
Then: 

\vskip 4pt

G0) The subdiagram of the diagram of $F$ consisting of those pairs 
$(x,X)$, $x\in F(X)$ with $X$ finite is cofinal.

\vskip 4pt 

Axioms on $\C$:

G1) $\C$ has final object $1$ and fiber products (notice that this
implies the existence of all finite limits).

G2) $\C$ has initial object $0$, [finite] coproducts and quotient of
objects by a finite group.

G3) $\C$ has strict epi-mono factorizations. That is, given any
$f:X \to Y$, there is a factorization $X \mapright{e} I \mapright{i} Y$
where $e$ is a strict epimorphism and $i$ is a monomorphism.
Furthermore it is assumed that $I$ is isomorphic to a direct summand of
$Y$. That is, there exist a subobject $J \hr Y$ and an isomorphism
$I \amalg J \cong Y$.

\vskip 4pt

Axioms on $F$:

G4) $F$ is left exact, that is, it preserves finite limits.

G5) $F$ preserves initial object, [finite] sums, quotients by
actions on finite groups and sends strict epimorphisms to surjections.

G6) $F$ reflects isomorphisms.

\vskip 6pt

We shall see that every object of $\C$ is a finite direct sum of
connected objects (see definition \equ(Dconnected)) and that the full
subcategory of non empty connected objects satisfy axioms C0) to C3) of
section III.

We start with an observation:

\proclaim{\Proposition (varios)}

1) The initial object $0$ is empty. That is, if we have an arrow
$X\to 0$, then $X=0$. We write $0=\emptyset$.

 Coproducts are disjoint and stable by pulling back.

2) Given any object $X \in \C$, $F(X) \cong \emptyset$ if and only if
$X\cong \emptyset$.

3) Given any finite objects $A,A_1,A_2$ and an isomorphism
$A\cong A_1 \amalg A_2$, if $A_1\cong A$ then $A_2\cong \emptyset$.

4) The functor $F$ preserves and reflects monomorphisms.

\endproclaim

\proclaim{Proof} \endproclaim

1) and 2) follows immediately by G5) and G6).

To see 3) we apply the functor $F$ and by G5) we get a bijection
of finite sets, $F(A) \cong F(A_1) \amalg F(A_2)$. 
It follows that $F(A_2) \cong \emptyset$. Thus by 2) $A_2\cong \emptyset$.

Finally, 4) follows from G4) and G6) since in any category an arrow
$u:X \to Y$ is a monomorphism if and only if the square

$$
\diagram
X \rto^{id} \dto_{id} & X \dto^{u}\\
X \rto_{u} &Y 
\enddiagram
$$
is a pullback.

\vskip 6pt

\proclaim{\Remark (preserves coprod)} Since pulling back along any 
arrow preserves coproducts, it follows immediately that if $A$ is 
connected, the representable functor $[A,-]$ preserves coproducts.
\endproclaim

\vskip 10pt

We consider now the diagram of the functor $F$. By G0) $F$ will also 
be the colimit of the subdiagram $\Tf$ consisting of the finite objects. 
More precisely, the objects of $\Tf$ are pairs $(x,X)$ where $x \in F(X)$
and $X$ is finite, and the arrows $(x,X) \mapright{f} (y,Y)$ are maps
$X\mapright{f} Y$ such that $F(f)(x)=y$.
Recall also that since $\C$ has finite limits and $F$ preserves
them, $\Tf$ is a filtered category. Thus $F$ is a pro-representable
functor and we denote $P$ the pro-object associated to $F$.
Therefore we have (by Yoneda) a correspondence

$$
\matrix
\underline{\q P \mapright{x} X \q}  \\
\underline{\q [X,-] \mapright{x} F \q} \\
\q x \in F(X) \q
\endmatrix
$$

\vskip 4pt

and the functor $F$ is ``representable'' by $P$. We write
$F(X) \cong [P,X]$.

\vskip 6pt

\proclaim{\Remark (P connected)} The fact that $F$ preserves [finite]
coproducts (G5) means in a sense that $P$ is connected.
\endproclaim

\proclaim {\Proposition (strict epi)} Every arrow $A \to B$ between
connected objects, with $A\not\cong \emptyset$ is an strict epimorphism.
\endproclaim

\proclaim{Proof} \endproclaim

Consider the factorization of $f$ given by G3):

$$
\diagram
A \drto|>>\tip \rrto^{f} && B \\
&I \urto|<<\tip   
\enddiagram
$$

Since $I$ is complemented and non empty, it follows that it's 
complement is empty. Thus $I \cong B$.

\proclaim{\Proposition (connected+strict epi)} Given any strict epimorphism $A \mapright{f} B$, 
if $A$ is connected so is $B$.
\endproclaim

\proclaim{Proof} \endproclaim

Let $B \cong B_1 \amalg B_2$.
Then by \equ(varios) 1) $A\cong  f^* B_1 \amalg f^* B_2$, where by $f^{*}$ we indicate 
pulling back along $f$.
It follows that $f^* B_1 \cong \emptyset$ or $f^* B_2 \cong \emptyset$.
Let $f^* B_1 \cong \emptyset$. We have then an strict epimorphism
$\emptyset \to B_1$. Since every arrow that starts in $\emptyset$ is a 
monomorphism (because $\emptyset$ is empty), it follows that $B_1 \cong 
\emptyset$.

\proclaim{\Definition (minimal)} An object of $\Tf$, $P\mapright{a} A$
is minimal if and only if it does not admit proper subobjects in
$\Tf$. That is, each time we have

$$
\diagram
P \dto_x \rto^{a} &A \\
X \urto_{u}|<<\tip  
\enddiagram
$$

with $u$ monomorphism, it follows that $u$ is an isomorphism,
thus $X \cong A$.
\endproclaim

\proclaim{\Proposition (connected=minimal)}
The following statements are equivalent:

i) $A$ is a finite non-empty connected object in the category $\C$.

ii) Every $P \mapright{a} A$ is minimal and $[P,A] \not\cong \emptyset$.

iii) There exist $P \mapright{a} A$ minimal.

\endproclaim

\proclaim{Proof} \endproclaim

i) $\Rightarrow$ ii)

Consider the following diagram

$$
\diagram
P \dto_x \rto^{a} &A \\
X \urto_{u}|<<\tip 
\enddiagram
$$

with $u$ monomorphism.

Since by G3) every subobject is complemented and $A$ is a non
empty connected object, it follows that $X \cong A$.

ii) $\Rightarrow$ iii) is clear.

iii) $\Rightarrow$ i)

Let $A\cong A_1 \amalg A_2$ and take $P \mapright{a} A$ minimal.
Since $F$ preserves coproducts, the arrow $a$ factors through $A_1$ or
$A_2$. Assume that it factors through $A_1$, thus $A \cong A_{1}$. Then by
proposition \equ(varios) item 3) it follows that $A_2\cong \emptyset$.

\vskip 6pt

We remark that to show i) $\Rightarrow$ ii) in this equivalence it is 
essential the whole strength of axiom G3), that is we need in the 
factorization that the subobject be complemented.

\proclaim{\Proposition (minimal cofinal)} The minimal objects
are cofinal in the diagram $\Tf$.
\endproclaim

\proclaim{Proof} \endproclaim

Let $P\mapright{x} X$ be an object of $\Tf$. If it is minimal we
are done. If not, we have a subobject $P \mapright{x_1} X_1$ such
that

$$
\diagram
&P \drto^{x} \dlto_{x_1} \\
X_1 \rrto|<<\tip  &&X 
\enddiagram
$$

If $P\mapright{x_1} X_1$ if minimal we are done. If not, we
have an $X_2$, etc...
It follows by proposition \equ(varios) item 4) and axiom G6) that a chain
like this stabilizes in a finite number of steps. We reach thus an
minimal object in $\Tf$ above $P\mapright{x} X$.

\proclaim{\Corollary (colimit)} The subdiagram $\Gf$ of $\Tf$ consisting 
of those objects $P\mapright{a} A$ with $A$ finite connected is 
cofinal in $\Tf$ and thus by G0) it is also cofinal in the whole 
diagram of $F$. Notice that this subdiagram is a poset by proposition 
\equ(strict epi), see remark \equ(poset).
\endproclaim

\proclaim{\Corollary (connected is finite)} Every connected object is 
finite.
\endproclaim

\proclaim{Proof} \endproclaim
Given $P \mapright{} X$, take 

$$
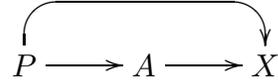
\diagram
P \rto \rrtou &A \rto &X
\enddiagram
$$

with $A$ finite connected. If $X$ is connected, then $A \to X$ is an 
strict epimorphism, thus $F(A) \to F(X)$ is surjective. It follows 
that $F(X)$ is finite because $F(A)$ is.

\proclaim{\Proposition (Gf meets)} The subdiagram $\Gf$ has meets
(finite infima).
\endproclaim

\proclaim{Proof} \endproclaim

Let $P\mapright{a} A$, $P\mapright{b} B$ and take $C$ connected
such that 

$$
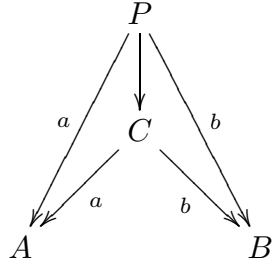
\diagram
& P \ddlto_a \ddrto^{b} \dto \\
& C \dlto^{a} \drto_{b} \\
A && B
\enddiagram
$$

By abuse of notation we indicate by the same label the
arrows from $C$. Take then the strict epi-mono factorization

$$
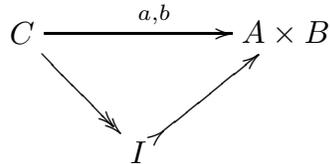
\diagram
C \drto|>>\tip \rrto^{a,b} && A \times B \\
&I \urto|<<\tip   
\enddiagram
$$

Since $C$ is connected, by proposition \equ(connected+strict epi), 
so is $I$. We shall see that
$P\to C\to I$ is the infima in $\Tf$ of $P\mapright{a} A$
and $P\mapright{b} B$. In fact, given

$$
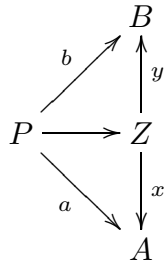
\diagram
& B \\
P \urto^{b} \rto \drto_{a} & Z \uto_{y} \dto^{x} \\
& A
\enddiagram
$$

with $Z$ connected, take the pullback $H$:

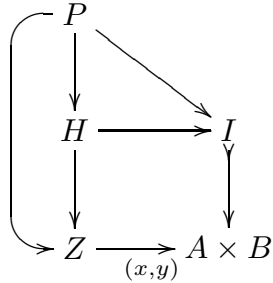
$$
\diagram
P \ddtol \dto \drto \\
H \rto \dto & I \dto|<<\tip \\ 
Z \rto_{(x,y)} & A \times B
\enddiagram
$$

\vskip 8pt

Since $P\to Z$ is minimal it follows that $H\cong Z$.
This finishes the proof.

\vskip 6pt

\proclaim{\Theorem (ConC)} The full subcategory of non empty
connected objects $\Con(\C) \hr \C$ together with the functor $F$
satisfies the axioms C0) to C3).
\endproclaim

\proclaim{Proof} \endproclaim

C0) By proposition \equ(strict epi).

C1) By axioms G2) and G5).

C2) By corollary \equ(connected is finite).

C3) By proposition \equ(Gf meets).

\vskip 10pt

We pass now to prove that every object in $\C$ is a [finite] coproduct of
connected objects.

\proclaim{\Proposition (coprod mono)} Given any two disjoint subobjects of $X$,
$A\to X$ and $B \to X$, the map $A \amalg B \to X$ is a monomorphism.
\endproclaim

\proclaim{Proof} \endproclaim

Apply the functor $F$ and use the corresponding result which holds
in $\Ens$. Use then the fact that $F$ reflects monomorphisms
(proposition \equ(varios) item 4).

\proclaim{\Proposition (coprod disjoint)} Given any two
connected subobjects $A\rightarrowtail X$ and $B \rightarrowtail X$, 
then $A\cap B\cong \emptyset$
or there exist an isomorphism $A \cong B$ such that

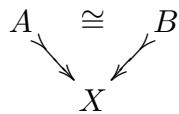
\spreaddiagramrows {-1pc}
\spreaddiagramcolumns {-1,3pc}
$$
\diagram
A \drto|<<\tip & \cong & B \dlto|<<\tip \\
&X
\enddiagram
$$
\spreaddiagramrows {1pc}
\spreaddiagramcolumns {1,3pc}

\endproclaim

\proclaim{Proof} \endproclaim

Consider the pullback diagram and the strict epi-mono factorization given
by G3) as indicated in the following diagram:

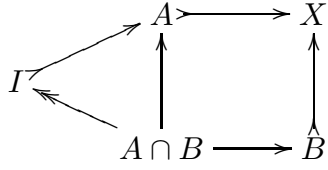
\spreaddiagramrows {-1,3pc}
$$
\diagram
&A \rto|<<\tip & X \\
I \urto|<<\tip \\
&A\cap B \uuto \rto \ulto|>>\tip & B \uuto|<<\tip  
\enddiagram
$$
\spreaddiagramrows {1,3pc}

Since $A$ is connected, $I\cong \emptyset$ or $I\cong A$.
In the first case, $A\cap B\cong \emptyset$ and in the second
the map $A\cap B \hr A$ is a monomorphism and a strict epimorphism,
thus is an isomorphism.

\vskip 6pt

\proclaim{\Theorem (coprod)} Every [finite] object $X \in \C$ is a
[finite] coproduct of connected objects of $\C$.
\endproclaim

\proclaim{Proof} \endproclaim

Assume $X\not\cong \emptyset$. Given $P \mapright{x} X$, let $P\mapright{a_x} 
A_x$ be a factorization:

$$
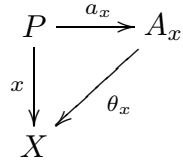
\diagram
P \dto_{x} \rto^{a_x} &A_x \dlto^{\theta_x} \\
X 
\enddiagram
$$

with $A_{x}$ connected. This factorization always exists by corollary
\equ(colimit).

Take the strict epi-mono factorization of each $\theta_x$ given by G3):

$$
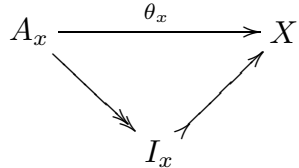
\diagram
A_x \drto|>>\tip \rrto^{\theta_x} && X \\
&I_x \urto|<<\tip  
\enddiagram
$$

By proposition \equ(connected+strict epi), the objects
$I_x$ are connected and by proposition \equ(coprod disjoint)
we can take a subfamily $I_l$
with $l \in J \subseteq [P,X]$ such that:

i) If $l\neq s$ then $I_l \cap I_s\cong \emptyset$

ii) $\forall$ $x\in [P,X]$, there exists $l \in J$ and an isomorphism
such that

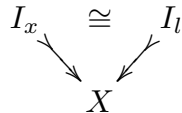
\spreaddiagramrows {-1pc}
\spreaddiagramcolumns {-1,3pc}
$$
\diagram
I_{x} \drto|<<\tip & \cong & I_{l} \dlto|<<\tip \\
&X
\enddiagram
$$
\spreaddiagramrows {1pc}
\spreaddiagramcolumns {1,3pc}

From i) and proposition (\equ(coprod mono)) it follows that the map
$\coprod_{l\in J} I_l \mapright{\lambda} X$ is a monomorphism,
thus $F(\lambda)$ is injective. We shall see now that $F(\lambda)$
is surjective. Then, by G6) it will follow that $\lambda$ is an
isomorphism.

In fact, let $x \in F(X)$, that is $P\mapright{x} X$. We have
$P \to I_x \to X$. It follows immediately from ii) above
that $x$ comes from some $l \in F(I_l)$,

$$
\diagram
P \dto_x \rto^l &I_l \dlto|<<\tip \\
X 
\enddiagram
$$

This finishes the proof. Clearly, if $X$ is finite, this coproduct is 
finite.

\proclaim{\Proposition(pi action)} Given any $X \in \C$, the set 
$[P,X]=F(X)$ has a continuous left action by the group $\pi=Aut(P)^{op}$
(which now is not transitive when $X$ is not connected).
\endproclaim

\proclaim{Proof} \endproclaim
Write $X\cong \coprod_{i\in I} X_{i}$ with $X_{i}$ connected. It 
follows from theorem \equ(ConC) (and proposition \equ(action)) that 
each $[P,X_{i}]$ has a continuous left action by the group $\pi$. The 
proof follows the immediately since the functor $[P,-]$ preserves 
coproducts. In fact, consider the following diagram:

$$
\diagram
\pi \times [P,X] \rdashed|>\tip & [P,X] \\
\pi \times [P,X_{i}] \uto \rto &[P,X_{i}] \uto \\
\enddiagram
$$

where the vertical arrows are coproduct diagrams in $\Ens$. The 
horizontal doted arrow exist by the universal property, and it is 
continuous since the left vertical arrows are also a coproduct diagram 
in the category of topological spaces because $\pi$ is compact 
Hausdorff.

\vskip 6pt

We are now in conditions to establish Grothendieck's Theorem.

Let $[P,X]_\pi$ be the $\pi$-set defined in proposition \equ(pi action).
Then: 

\proclaim{\Theorem (Grothendieck)} Let $\C$ be a category such that 
axioms \equ(axiomasG) hold. Consider the functor:

\spreaddiagramcolumns {1pc}
$$
\diagram
\C \rto^{[P,-]_\pi} & \Ens^{\pi} \\
\enddiagram
$$
\spreaddiagramcolumns {-1pc}

defined in proposition \equ(pi action). Then this functor is an 
equivalence of categories.

Let $f\C$ be the full subcategoy of finite objects. Then the functor 
$[P,-]_{\pi}$ restricts into the functor:

\spreaddiagramcolumns {1pc}
$$
\diagram
f\C \rto^{[P,-]_\pi} & f \Ens^{\pi} \\
\enddiagram
$$
\spreaddiagramcolumns {-1pc}

which also establishes an equivalence of categories.
\endproclaim

\proclaim{Proof} \endproclaim

We have commutative diagrams of categories and functors:

$$
\diagram
\C \rto^{[P,-]_\pi}& c \Ens^\pi &  f\C \rto^{[P,-]_\pi}& f c \Ens^\pi\\
\Con(\C) \uto|<\ahook \rto^{[P,-]_\pi} & ct \Ens^\pi \uto|<\ahook  &
\Con(\C) \uto|<\ahook \rto^{[P,-]_\pi} & ct \Ens^\pi \uto|<\ahook 
\enddiagram
\Eq(equiv)
$$

where by $\Con(\C)$ we indicate the full subcategory of non empty connected 
objects.

It is immediate to see that every object in $c \Ens^{\pi}$ (resp. 
$fc\Ens^{\pi}$) is a coproduct (respectively finite coproducts) of 
transitive actions. Define the functors:

\spreaddiagramcolumns {1pc}
$$
\diagram
\C & c \Ens^\pi \lto^{P\times_\pi (-)} & \C & 
c \Ens^\pi_{< \infty} \lto^{P\times_\pi (-)} \\
\enddiagram
$$
\spreaddiagramcolumns {-1pc}

using coproducts in $\C$ (resp. finite coproducts in $f\C$) 
and the fact that it is already defined in $ct \Ens^{\pi}$.
Since the functors:

\spreaddiagramcolumns {1pc}
$$
\diagram
\C \rto^{[P,-]_\pi} & c \Ens^{\pi} &
\C \rto^{[P,-]_\pi} & c \Ens_{< \infty} ^{\pi} \\
\enddiagram
$$
\spreaddiagramcolumns {-1pc}

preserve coproducts (resp. finite coproducts), the proof follows from
theorem \equ(coprod) and the fact that the arrows below in the 
diagrams \equ(equiv) are an equivalence of categories by theorem \equ(ConC)
and theorem \equ(adjoint2) in section III.

\vskip 10pt

\subheading{V - All actions of a discrete group}
\numsec=5\numtheo=1\numfor=1

\vskip 10pt

In section II we developed the representable and connected case, which 
characterize the category of transitive actions of a discrete group. We 
restrict to the connected case for two reasons. First, it develops in an 
abstract setting exactly the same techniques that are utilized by 
topologists in the theory of covering spaces and its relation with the 
fundamental group. This clearly shows where the ideas come from. Second, it 
is the only case needed for Grothendieck's development of the general 
case (non representable non connected), as we show in sections III and IV. Here 
quotients by group actions on connected objects play a central role, and 
coproducts are only used to extend the constructions to the non connected 
objects. 

However, it seems to us that in the representable case different 
phenomena are also behind the result. Concretely, a direct 
characterization of the category (now a topos) of all actions of a discrete 
group can be obtained, rather than derive it from the connected case (as we 
do for the non representable case in section IV). The proof uses some 
general categorical techniques that have its own interest and apply to many 
other situations (like Beck's tripability theorem, Giraud theorem on 
characterization of topoi, and in the additive case, representation theorems 
for abelian categories and the Morita equivalences).

\vskip 6pt
 
Let $\C$ be any category and $A \in \C$ an object. Let
$G=[A,A]^{op}$ the opposite monoid of the monoid of endomorphisms
of $A$.

\proclaim{\Proposition (action3)}
The monoid $G$ acts on the left on the sets $[A,X]$
of morphisms in $\C$ from $A$ to any other object $X \in \C$:
\endproclaim

$$
\matrix
G&\times&[A,X]&\longrightarrow&[A,X] \\
g&&x&\longrightarrow&gx=x\circ g
\endmatrix
$$

where ``$\circ$'' is composition of arrows in $\C$.

\vskip 4pt

We have in this way a functor:

\spreaddiagramcolumns {1pc}
$$ 
\diagram
\C \rto^{[A,-]_G} & \Ens^{G} 
\enddiagram
$$
\spreaddiagramcolumns {-1pc}

where by $[A,X]_G$ we indicate the set $[A,X]$ with this
action of $G$ and $\Ens^G$ indicates the category of $G$-sets.

\vskip 4pt

We assume now that the category $\C$ has coproducts and coequalizers. It 
is well known and easy to see that the representable functor $[A,-]$ has a 
left adjoint (the ``tensor'' in closed category terminology) that we denote here 
$A\bullet (-)$, and that it is given by the formula 

$$A\bullet S = \coprod_S A$$
 
Clearly, by definition of coproducts, there is a bijection:

$$
\matrix
\underline{\q A\bullet S \to X \q} & \text{natural in $X$} \\
\q S \to [A,X] \q
\endmatrix
\Eq(bijection2)
$$

In the case that $\C = \Ens^G$, for any monoid $G$, and $A = G$ with its 
canonical action on itself, the representable functor $[G,-]$ is just the 
underline set functor
(G is the free action on one generator), $[G,E] = |G|$,  and the 
functor $G\bullet (-)$ is conveniently given by the formula:

$$G\bullet S = G \times S \q \text{with the action} \q g(f, s) = (gf, s)$$

We have then the following diagram of categories and functors:

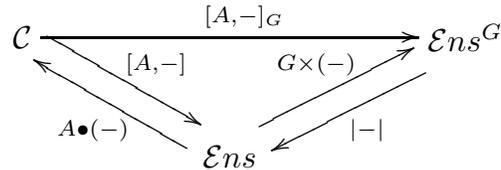
\spreaddiagramcolumns {2pc}
$$
\diagram
\C \drto<1ex>^{[A,-]} \rrto^{[A,-]_{G}} && \Ens^{G} \dlto<1ex>^{|-|}  \\ 
&\Ens \ulto<1ex>^{A \bullet (-)} \urto<1ex>^{G\times (-)}
\enddiagram
$$
\spreaddiagramcolumns {-2pc}

where the triangle formed by the right adjoints commutes.

In this situation, if the category $\C$ has coequalizers, there is a left 
adjoint $A \times_G (-)$ for the functor $[A,-]_G$ that it is given by an 
specific construction in terms of the other two adjointness in the triangle 
(see [3]).
In the particular situation considered here we have:

\proclaim{\Proposition (leftadjointV1)}
Let $\C$ be a category with coproducts and coequalizers (thus, all colimits), 
and let $A$ in $\C$. Then, the functor $[A, -]_G$ has a left 
adjoint $A \times_G (-)$ defined, given $E$ in $\Ens^{G}$, as the following 
collective coequalizer:

$$
\diagram
A \drto^{\lambda_{x}} \\
& A\bullet |E| \rto^{q} & A\times_{G} E \\
A \uuto_{g} \urto_{\lambda_{gx}} \\
\enddiagram
$$

That is, the arrow $A\bullet |E| \mapright{q} A \times_G E$ is universal with 
respect to the equations  $q \circ \lambda_x \circ g = q \circ \lambda_{gx}$, 
one for each $x \in E$ and $g \in [A,A]$. 
\endproclaim

\proclaim{Proof} \endproclaim

The proof is immediate, just check that in the bijection \equ(bijection2), 
an arrow $A\bullet |E| \mapright{f} X$ factors through $A \times_G E$ 
(that is, it satisfies all equations 
$f \circ \lambda_x \circ g = f \circ \lambda_{gx}$) if and only if the 
corresponding arrow $E \to [A, X]_G$ is a morphism  of actions.

\vskip 6pt

We set now the following definition :

\proclaim{\Definition (Dgenerator)} An object $A$ in a category with 
coproducts $\C$ is called a generator if for every $X \in \C$ 
the following diagram:

$$
\diagram
A \drto^{\lambda_{x}} \\
& A\bullet [A,X] \rto^>>>>>{\varepsilon} & X \\
A \uuto_{g} \urto_{\lambda_{gx}} \\
\enddiagram
$$

\vskip 3pt

where $\varepsilon$ is the unique map such that $\varepsilon \circ 
\lambda_x=x$, is a collective coequalizer. 
In particular, every object is a quotient of a sum 
of copies of $A$. 
\endproclaim

Notice that the map $A\bullet [A, X] \mapright{\varepsilon} X$ induces 
precisely the counit of the adjunction $A \times_G [A, X]  \to  X$. 
Thus, $A$ is a generator if and only if this counit is an isomorphism. 
This notion is some times called ``dense generator''.

\proclaim{\Proposition (Pgenerator)}
Let $\C$ be a category with coproducts and coequalizers and $A \in \C$. 
Then:

\item{i)} If $A$ is a generator then the representable functor $[A, -]$ 
reflects isomorphisms.
\item{ii)} Assume the representable functor $[A, -]$ preserves coproducts 
and coequalizers. Then, if it reflects isomorphisms, $A$ is a generator.
\endproclaim

\proclaim{Proof} \endproclaim

\item{i)} Let $X \mapright{f} Y$ and consider the commutative square:

$$
\diagram
A\bullet [A, X] \dto_{A\bullet f_*} \rto & X \dto_{f} \\
A\bullet [A, Y]  \rto & Y
\enddiagram
$$

where the horizontal arrows are collective coequalizers with respect to the 
equations in definition \equ(Dgenerator). If $f_* = [A, f]$ is a bijection, 
$A\bullet f_*$ is an isomorphism, and one can readily check that the bijection 
$f_*$ establishes also a correspondence between the respective equations.
  
\item{ii)} Let $A\bullet [A, X] \mapright{q} Q$ be the collective coequalizer, 
and consider:

$$
\diagram
A\bullet [A, X] \rrtou^{\varepsilon} \rto^>>>>>{q} &Q \rto^{\eta} &X 
\enddiagram
$$

We shall see that $\eta$ is an isomorphism. Applying the functor $[A, -]$ we 
have:

$$
\diagram
[A, A\bullet [A, X]] \rrtou^{\varepsilon_{*}} \rto^>>>>>{q_{*}} &[A,Q] 
\rto^{\eta_{*}} &[A,X] 
\enddiagram
$$

Where $q_*$ is the respective coequalizer in the category of sets. Since for 
every $x \in [A, X]$ we have $\varepsilon_*(\lambda_x) = x$, we have 
that $\eta_*$ is surjective. To see that it is also injective we have to see 
that $\varepsilon_*$ does not identify more than $q_*$.
Let $u, v \in [A, A\bullet [A, X]]$ be such that $\varepsilon \circ u = 
\varepsilon \circ v$. Since $[A, -]$ preserves coproducts, 
$u = \lambda_x \circ h$, $v = \lambda_y \circ g$, for some $x, y \in [A, X]$ 
and $h, g \in G$. 

In the next chain of equations, each one follows immediately from the 
preceding one:

$\varepsilon \circ u=\varepsilon \circ v \q \Rightarrow \q 
\varepsilon \circ \lambda_{x} \circ h=\varepsilon \circ \lambda_{y} 
\circ g \q \Rightarrow \q 
\varepsilon \circ \lambda_{hx} =\varepsilon \circ \lambda_{gy} \q \Rightarrow \q 
hx =gy \q \Rightarrow \q 
\lambda_{hx} =\lambda_{gy} \q \Rightarrow \q 
q \circ \lambda_{hx} =q \circ \lambda_{gy} \q \Rightarrow \q 
q \circ \lambda_{x} \circ h =q \circ \lambda_{y} \circ g \q \Rightarrow \q 
q \circ u=q \circ v$.


\proclaim{\Nada (axiomasR)} {\bf Axioms for the representable case}
\endproclaim

Consider a category $\C$ and an object $A \in C$.

Axioms on $\C$:

R1) $\C$ has a terminal object and pullbacks (thus all finite limits).

R2) $\C$ has coequalizers.

R3) $\C$ has coproducts.

\vskip 4pt

Axioms on $A$ (in terms of the representable functor $[A, -]$):

R4) $[A, -]$ preserves coequalizers. 

R5) $[A, -]$ preserves coproducts.

R6) $[A, -]$ reflects isomorphisms.

\proclaim{\Theorem (teoremaR)} Let $\C$ be any category and $A \in \C$ as in 
\equ(axiomasR). Then the left adjoint $A \times_G (-) \dashv [A, -]_G$
exists and establishes an equivalence of categories.
\endproclaim

\proclaim{Proof} \endproclaim
By R6 and proposition \equ(Pgenerator) the counit of the adjunction in $\C$
is an isomorphisms. Since $G$ is itself a generator (in the sense of definition 
\equ(Dgenerator)) for the category $\Ens^G$, and by R4 and R5 the monad 
(composite of the two functors) in $\Ens^G$ preserves coproducts and 
coequalizers, it is enough to check that the unit is an isomorphisms in the 
object $G$. But this is clear. In fact, using again 
that $A$ is a generator, we have $A \times_G G$ is isomorphic to $A$, thus 
$[A, A \times_G G]$ is isomorphic to $[A, A] = G$.       

\vskip 6pt

It is interesting to observe that it is not necessary to assume any 
exactness properties in the category $\C$. However, if we want to write the 
axioms R4 and R5 as properties of the object $A$ (projective and connected), 
then we need some of the exactness properties characteristic of a topos (or 
pretopos).

We recall (in \equ(strictepi), \equ(epimonofact), \equ(DPregular), and
\equ(excoequalizers)) from [[1], I, 10] the definition of strict 
epimorphism and some of its properties. The interested reader can 
also consult [2].

\proclaim{\Defprop (strictepi)}

1) An arrow $X \mapright{f} Y$ in any category $\C$ is an strict 
epimorphism if given any compatible arrow $X \mapright{g} Z$, there 
exists a unique $Y \mapright{h} Z$ such that $g = h \circ f$. $g$ is 
compatible if for all 
$\C  \mapright{x} X$, $\C  \mapright{y} X$
with $g \circ x = g \circ y$, then also $f \circ x = f \circ y$.

2) It immediately follows that:

           $$\text{strict epimorphism  +  monomorphism   =   isomorphism}$$

3) The kernel pair of an arrow  $X \mapright{f} Y$ is the pull-back of  
$f$  with itself, and it will be denoted $R_f$. The subobject 
$R_f \subset X \times X$ is always an equivalence relation. 

4) An equivalence relation is called effective if the quotient exists and it 
becomes the kernel pair of this quotient.   

5) An strict epimorphism with a kernel pair is called effective. 


6) Given any arrow $X \mapright{f} Y$, the diagram 
\spreaddiagramcolumns {1pc}
$$
\diagram
R_{f} \rto<1ex> \rto<-1ex>& X \rto^{f} &Y \\
\enddiagram
$$
\spreaddiagramcolumns {-1pc}

is a coequalizer if and only if $f$ is an strict epimorphism.

7) A diagram  
\spreaddiagramcolumns {1pc}
$$
\diagram
Z \rto<1ex>^{g} \rto<-1ex>_{h}& X \rto^{f} &Y  \\
\enddiagram
$$
\spreaddiagramcolumns {-1pc}

is a coequalizer if and only if 
$f$ is an strict epimorphism, and $R_{f}$ is the smallest effective equivalence 
relation that factories the arrow $(g, h):Z \mapright{} X \times X$.   
\endproclaim

\proclaim{\Proposition (epimonofact)}
(strict epi - mono factorizations, Tierney-Kelly).
Let $\C$ be any category with pull-backs and quotients of equivalence relations 
which are stable under pulling-back (equivalently, strict epimorphisms are 
stable under pulling back). Then $\C$ has strict epi-mono factorizations. 
That is, given any $f:X \to Y$, there is a factorization 
$X \mapright{e} I \mapright{i} Y$
where $e$ is a strict epimorphism and $i$ is a monomorphism.
\endproclaim

\proclaim{Proof} \endproclaim
Let $X \to Z$ be the quotient of $X$ by the kernel pair of $f$. Clearly there 
is a factorization $X \to Z \to Y$. From the stability of strict epimorphisms 
it readily follows that the arrow $Z \to Y$ is a monomorphism.

\proclaim{\Defprop (DPregular)}

1) A category $\C$ is Regular if it has finite limits and strict epi - mono 
factorizations (thus, direct images) stable under pull-backs. 

2) A relation in $X$ is a subobject of $X \times X$. Using pull-backs 
and images the composition of relations is defined in the usual way as 
in the category of sets. 

3) If $\C$ has unions of denumerable chains of subobjects which are stable 
under pulling-back, the transitive closure of any relation is constructed as 
the union of iterated compositions. Thus, if in addition there are finite unions
of subobjects (to make a relation reflexive and symmetric), any relation is 
contained in a smallest equivalence relation, which is generated as described 
above. Notice that denumerable stable coproducts are sufficient, since unions 
are the direct image of the induced map from the coproduct. 
\endproclaim

\proclaim{\Proposition (excoequalizers)} Let $\C$ be a Regular Category such 
that every 
equivalence relation has a quotient and every denumerable chain of subobjects 
has an stable union (or denumerable stable coproducts exists in $\C$). 
Then, $\C$ has coequalizers. 
If in addition, equivalence relations are effective, then any finite limit 
preserving functor $F: \C \to \E$ into any other such category will preserve 
coequalizers if it preserves strict epimorphisms and stable 
denumerable unions (or preserves stable denumerable coproducts).     
\endproclaim

\proclaim{Proof} \endproclaim 
Given any pair of arrows $Z \mapright{g} X$, $Z \mapright{h} X$,
let $R(g, h)$ be the 
equivalence relation generated (as in \equ(DPregular).3) by the image 
of $X \mapright{(g, h)} X \times X$, and $X \mapright {f} Y$ its quotient. 
If an arrow $X \mapright{s} H$ coequalizes $g$ and $h$, $(g, h)$ 
factors through its kernel pair $R_{s} \subset X \times X$. Thus 
$R(g,h) \subset R_{s}$, which implies there is a map $Y \to H$ showing that 
$X \mapright{f} Y$ is the coequalizer of $g$ and $h$. For the second 
assumption, notice that if equivalence relations 
are effective, 7) in \equ(strictepi) says that $R_{f}$ is just the 
smallest equivalent relation that factories $(g, h)$. Then, the proof is clear. 

\vskip 10pt

\proclaim{\Nada (axiomasE)} {\bf Axioms for the representable case 
(second version)}
\endproclaim

Consider a category $\C$ and an object $A \in C$.

Axioms on $\C$   

E1) $\C$ has a terminal object and pullbacks.

E2) $\C$ has quotients of equivalent relations which are then the kernel 
pair of this quotient, and this quotient is stable by pulling back (i.e.
equivalence relations are effective and universal).

E3) $\C$ has an initial object and coproducts stable by pulling back.

Notice that we do not require that coproducts be disjoint.

Axioms on $A$:

E4) $A$ is projective (see \equ(Dproyective)).

E5) $A$ is connected (see \equ(Dconnected)).

E6) $A$ is a generator (see \equ(Dgenerator)).

Notice that E6) implies that $A$ is non empty.

\vskip 8pt

\proclaim{\Defprop (Dconnected)}
An object $A$ in a category $\C$ is connected if and only if
it satisfies any of the two following equivalent statements:

i) If $A \cong A_1 \amalg A_2$, then $A_1\cong \emptyset$
or $A_2 \cong \emptyset$.

ii) If $A \cong \coprod_{i\in I} A_{i}$, with $I$ any set, then there 
exists an $i \in I$ such that $A \cong A_{i}$ and $A_{j} \cong \emptyset$ 
for every $j\neq i$.

\endproclaim

\proclaim{Proof} \endproclaim

We can assume $A \not\cong \emptyset$. For each $i$ we have $A \cong A_{i} 
\amalg \coprod_{j\neq i} A_{j}$. Clearly, there exists $i$ such 
that $A_{i} \not\cong \emptyset$. Then 
$\coprod_{j\neq i} A_{j} \cong \emptyset$. Thus, $A_{j}\cong 
\emptyset$ for all $j$ (consider the arrow $A_{j} 
\mapright{\lambda_{j}} \coprod_{j} A_{j}$).

\vskip 4pt

\proclaim{\Remark (preservescoprod2)} By axiom E3) pulling back along any 
arrow preserves coproducts, then it follows immediately that if $A$ is 
connected, the representable functor $[A,-]$ preserves coproducts.
If the initial object is empty and coproducts are disjoint, then the other 
implication also holds.
\endproclaim

\proclaim{\Definition (Dproyective)} 
By projective we mean the assumption that the functor 
$[A, -]$ preserves strict epimorphisms (this is weaker that the usual 
meaning of projective, which is defined with respect to all epimorphisms).
\endproclaim

\proclaim{\Theorem (teoremaE)} Let $\C$ be any category and $A \in \C$ as in 
\equ(axiomasE). 
Then the left adjoint $A \times_{G} (-) \dashv [A, -]_{G}$ exists and 
establishes an equivalence of categories.
\endproclaim

\proclaim{Proof} \endproclaim
We shall see that the axioms in theorem \equ(teoremaR) hold. R1) = E1), and from 
propositions \equ(epimonofact) and \equ(excoequalizers) we have R2). 
E4) and E5) 
imply R4) by remark \equ(preservescoprod2) 
and proposition \equ(excoequalizers), 
and E5) gives R5) also by remark \equ(preservescoprod2)
Finally, E6) implies R6 by proposition \equ(Pgenerator).
   
\vskip 11pt

{\bf Group actions}

\vskip 6pt

Group actions come into the scenery when the monoid $[A,A]$ is actually a 
group. That is, every endomorphism of $A$ is an isomorphism, $Aut(A) = [A, A]$.
In this case, the theorems above are much simpler. Notably, there is no 
need to develop all the calculus of relations in regular 
categories, neither to consider equivalence relations, transitive 
closures and coequalizers in general. In 
place of all this, quotients by group actions do all the work.

\vskip 6pt

To guide the reasoning it is useful to consider now the elementwise 
description of coproducts which consist on writing informally:

$$
A\bullet S = \{(a, x) / x \in S \q \text{and} \q 
a \in \text{the x-copy of $A$}\}
$$

Then the object $A \times_G E$, $G = Aut(A)^{op}$, 
can be described informally as the object $A\bullet |E|$ divided by 
the equivalence relation (clearly reminiscent of the tensor 
product construction) defined by pairs $(ag, x) = (a, gx)$ ($G$ acts on the 
right on $''a''$ because it is the opposite group $Aut(A)^{op}$). 
This equivalence relation is thus $(b, x) = (a, y)$ if and only if 
there exists $g \in G$ such that $b = ag$ and $y = gx$.
There is a left action of $G$ in $A\bullet |E|$ given by 
$g(a, x) = (g^{-1}a,gx)$, and the orbits of this action 
define precisely the equivalence classes of the relation above.

It is instructive to see directly the consistency of this construction with 
the one given for transitive actions in section II (both constructions must
coincide by the universal property of the adjoint). Given an $x_{0} \in E$, 
we have the equations 
$q \circ \lambda_{x_{0}} \circ h = q \circ \lambda_{x_{0}}$, for all 
$h \in H = Fix(x_0)$. Thus the map $\lambda_{x_{0}}$
induces an injective map $A/H \to A \times_G E$. 
When the action is transitive, it 
is immediate to see using the description above that it is also surjective.    
   
Formally:

\vskip 4pt

\proclaim{\Proposition (GA)}
There is a left action of $G$ on $A\bullet |E|$ defined by the following 
diagrams (one for each $g \in G$):

$$
\diagram
A\bullet |E| \rto & A\bullet |E| \\
A \uto^{\lambda_{x}} \rto^{g^{-1}} & A \uto_{\lambda_{x}}
\enddiagram
$$

and the quotient by this action (definition \equ(quotient)) is precisely 
the collective coequalizer that defines  $A \times_G (-)$.
\endproclaim

\vskip 4pt

\proclaim{\Theorem (teoremaGR)}
Let $\C$ be any category and $A \in \C$ as in \equ(axiomasR) but with 
axioms R2) and R4) respectively replaced by:

R'2) $\C$ has quotients of objects by groups of automorphisms.           

R'4) $[A, -]$ preserves quotients of objects by groups of automorphisms.

Assume that $[A,A]=Aut(A)$, then the left adjoint 
$A \times_G (-) \dashv [A, -]_G$
exists and establishes an equivalence of categories.
\endproclaim

\proclaim{Proof} \endproclaim

It follows the same lines that the proof of theorem \equ(teoremaR). 
By the proposition 
above and R'2 we can construct the functor $A \times_G (-)$. 
From definition \equ(Dgenerator) 
it is clear that when $G$ is a group R'4 and R5 imply that it is enough to check
that the unit is an isomorphisms in the object $G$, which we do as in 
theorem \equ(teoremaR).

\vskip 4pt

\proclaim{\Defprop (GDP)}
Let $\C$ be any category with finite limits and coproducts
and let $H$ be a group that acts by automorphisms in an object $A \in C$. 
Then the quotient $A \mapright{q} A/H$ (see definition \equ(quotient))
fits into a coequalizer:

\spreaddiagramcolumns {1pc}
$$
\diagram
A \bullet H  \rto<1ex>^{\nabla} \rto<-1ex>_{\nabla_H} 
& A \rto^{q} & A/H\\
\enddiagram
$$
\spreaddiagramcolumns {-1pc}

where $\nabla$ is the codiagonal and $\nabla_H$ is defined by:

$$
\diagram
A \bullet H \rto^{\nabla_{H}} & A \\
A \urto_{h} \uto^{\lambda_{h}}
\enddiagram
$$

\vskip 4pt

We say that the action of $H$ on $A$ is effective if the quotient exists and 
the map 

$$
\diagram
A \bullet H \rrto^{(\nabla,\nabla_{H})} && R_{q} 
\subset A \times A
\enddiagram
$$
 
is an strict epimorphism.
\endproclaim

\vskip 4pt

\proclaim{\Proposition (Gexcoequializers)}
Let $\C$ be any category with finite limits, coproducts and 
where actions of groups by automorphisms are effective. Then any finite limit 
preserving functor $F:\C \to \E$ into any other such category will preserve 
quotients by actions of groups if it preserves strict epimorphisms and 
coproducts.
\endproclaim

\proclaim{Proof} \endproclaim

Consider the diagram

\spreaddiagramcolumns {1pc}
$$
\diagram
& R_{q} \dto<1ex> \dto<-1ex> \\
A \bullet H \urto|>>\tip^{e} \rto<1ex>^{\nabla} 
\rto<-1ex>_{\nabla_H}& A \rto|>>\tip^{q} & A/H\\
\enddiagram
$$
\spreaddiagramcolumns {-1pc}                       
    
By the assumptions made this diagram goes by $F$ into a diagram:

\spreaddiagramcolumns {1pc}
$$
\diagram
& R_{Fq} \dto<1ex> \dto<-1ex> \\
FA \bullet H \urto|>>\tip^{Fe} \rto<1ex>^{\nabla} 
\rto<-1ex>_{\nabla_H}& FA \rto|>>\tip^{Fq} & F(A/H)\\
\enddiagram
$$
\spreaddiagramcolumns {-1pc}                       

That $F(q)$ is the coequalizer of $\nabla$ and $\nabla_{H}$ 
follows immediately from 
the fact that $F(e)$ and $F(q)$ are both strict epimorphisms.

\vskip 4pt

\proclaim{\Theorem (teoremaGE)}
Let $\C$ be any category and $A \in \C$ as in \equ(axiomasE) but with
E2) replaced by:

E'2) In $\C$ actions of groups by automorphisms are effective. 

Assume that $[A, A] = Aut(A)$, then the left adjoint 
$A \times_{G} (-) \dashv [A, -]_{G}$ exists and 
establishes an equivalence of categories.
\endproclaim

\proclaim{Proof} \endproclaim 
Construct the functor $A\times_{G} (-)$ as in theorem \equ(teoremaGR).
To finish the proof it remains to see that $[A,-]$ preserves 
quotients by group actions. But this is done in the proposition above.

\vskip 10pt

\subheading{VI - Final Comments}
\numsec=6\numtheo=1\numfor=1

\vskip 10pt

We have seen how Grothendieck's interpretation of the Galois's 
Galois Theory (and its natural relation with the fundamental 
group and the theory of coverings) allows the development of a 
categorical theory much beyond the original Galois's scope.
Here are some further results essentially proved in the previous 
sections and some comments about them.

\vskip 6pt

\proclaim{\Nada (IIyV)}{\bf On sections II and V}
\endproclaim

\vskip 4pt

We can easily generalize the results in section V by replacing the 
single object $A$ by a (small) set of objects and its associated 
full subcategory  $\A  \subset \C$.
In this case the category $\Ens^{G}$ becomes the presheaf category
$\Ens^{\A^{op}}$ and the restriction of the Yoneda embedding 
$h: \C \to  \Ens^{\A^{op}}$ corresponds to the functor $[A,-]_G$. 
This generalizes proposition \equ(action3).
The left adjoint to h is given by a Kan extension 
that we denote $k:\Ens^{\A^{op}} \to \C$. Given $E \in \Ens^{\A^{op}}$, 
consider its canonical diagram $\GE$ so that  

$$E=\dsize\colim_{(a,A) \in \GE}{[-, A]}$$
                            
Then $k(E)$ is the corresponding colimit taken in $\C$         

$$k(E)=\dsize\colim_{(a,A) \in \GE}{A}$$
                                   
The adjointness is easily seen from this definition. This generalizes 
Proposition \equ(leftadjointV1). 

\vskip 4pt

We pass now to the axiom set \equ(axiomasE). Axioms E1), E2) and E3) are 
respectively conditions a), c) and b) in Giraud's Theorem [[1], IV, 1.2]. 
Axiom E6) says that the category $\C$ has an (small) generating set, 
which is condition d) in that theorem. We have in addition axioms E4) and 
E5), 
which say that the objects in the generating set are projective and 
(non-empty) connected. Thus, theorem \equ(teoremaE) is Giraud's Theorem 
in the particular case in which E4) and E)5 hold. It gives a characterization 
of categories of presheaves. Compare also with [[1], IV, 7.6], where condition 
[ii): the full subcategory of (non empty) connected projective objects 
is generating] in exercise 7.6 d) consists precisely on the axioms E4), 
E5) and E6) together. This characterization is due to J. E. Ross.
 
Giraud's Theorem can be analyzed as follows: 
Assume axioms E1), E2), E3) and E6) (not E4) and E5)), that is, Giraud's 
assumptions. Then we have an adjunction as indicated above 
(proposition \equ(leftadjointV1). Axiom E6) says (tautologically) that 
the counit of this adjunction is an isomorphism, making $\C$ a full subcategory
 of 
$\Ens^{\A^{op}}$. How can we make the unit of the adjunction in $\Ens^{\A^{op}}$
be an isomorphism 
without E4) and E5) ?, (as we do in theorem \equ(axiomasE)). 
Giraud's answer is to force this unit to be an isomorphism by means 
of the appropriate Grothendieck topology (here it is probably needed 
that coproducts be disjoint in addition to stable).

All this relativices to the additive case, where, in the place of 
sets, the category of abelian groups is the base category. In fact, it was this
case that was developed first. The statements and their proofs are 
essentially the same. A conspicuous difference is that now preservation
of coproducts by the representable functor 
$[A, -]$ means that the object $A$ is ``small'', instead of ``connected''. 
A coproduct of two small objects is still small, and this implies 
that the category of modules over a ring does not characterize the 
ring  (while the category of actions of a monoid does characterize the 
monoid). An earlier instance of the additive case of theorem \equ(teoremaE)
 are the Morita 
equivalences. The abstract additive theorem \equ(teoremaE) was 
known I imagine by the end of the fifties. It is explicitly stated 
for example in [[5], Ch 4, F]. The generalization to several objects 
discussed above in [[5], Ch 5, H]. 

\vskip 10pt
       
\proclaim{\Nada (IIIyIV)} {\bf On sections III and IV}
\endproclaim

\vskip 4pt

In section III, theorem \equ(adjoint2), we proved a characterization 
of the category of transitive actions of a profinite group by 
``passing into the limit'' (theorem \equ(adjoint1) in section II on 
transitive actions of a discrete group) in a filtered colimit of categories.
Using this, we show in section IV Grothendieck's theorem in [6]. 
This theorem (\equ(Grothendieck) in section IV) implies that a connected 
locally connected [[1], IV, 8.7, l)] boolean topos with a ``profinite point'' 
(that is, its inverse image functor satisfies axiom G0 in 
\equ(axiomasG)) is the topos of continuous $\pi$-sets for a 
profinite group $\pi$ (locally connected + boolean = atomic, 
see [[8], Ch. 8]). A careful analysis shows, in the light of axioms 
\equ(axiomasC) and proposition \equ(PC0), that axiom G6) can be omitted 
in \equ(axiomasG).

Our proof shows that Grothendieck' theorem is the result of 
``passing into the limit'' (theorem \equ(teoremaGR) in section V on all 
actions of a discrete group) in a filtered inverse limit of topoi. 
Warning: such a limit is not a colimit of categories and inverse 
image functors, but, nevertheless, we have such a colimit on the sites, 
(see [[1], VI, 8.2.3, 8.2.9], also [10]).

This theorem can be generalized to the case in which the point is not 
any more ``profinite''. Here it is necessary to consider the diagram of discrete 
groups without passing to its inverse limit group, since the projections 
from this limit can not any more be guaranteed to be surjective 
[[1], IV, 2.7]. Thus, the topological group is replaced by a whole system 
of discrete groups, called a ``progroup'' by Grothendieck. The 
corresponding theorem is a theorem of characterization of the category 
of all actions (as defined by Grothendieck) of a progroup. In this 
more general context the existence 
of the Galois Closure can not 
be proved as in proposition \equ(closure). It is necessary to introduce 
 an  stronger (than just finite 
limit) preservation property in the functor F. This condition 
(preservation of cotensors = arbitrary products of a same object) 
replace the finiteness condition G0. We will show how to do all this elsewhere.

In [[9], VIII, 3. theorem 1] the finiteness condition in theorem 
\equ(Grothendieck) it is also removed. This is done following a different 
line that the one sketched above. As it was stressed by Joyal since the
early seventies, topological spaces are replaced by ``spaces'', which are 
the dual objects of locales (sup-complete distributive lattices). 
The profinite (topological) group becomes an spatial group, which in 
general will not have sufficiently many points, and thus it will not be 
a topological group in the classical sense. An spatial group is a more 
general concept that a progroup. The topos of continuous $\pi$-sets of 
an spatial group satisfy the condition on preservation of cotensors 
presicely when the spatial group corresponds to a progroup.

Finally, let us point out that all this should be generalized taking 
a category (grupoid) of appropriate points, instead of a single point, as 
it is already developed in [6], and continued in [7]. It is not clear to us 
the relation of this with the principal theorem in [9] concerning continuous 
actions of an spatial grupoid.

\vskip 5pt
\vskip 5pt
\vskip 5pt
\vskip 5pt
\vskip 5pt
\vskip 5pt
\vskip 5pt
\vskip 5pt
\vskip 5pt
\vskip 5pt
\vskip 5pt
\vskip 5pt
\vskip 5pt
\vskip 5pt
\vskip 5pt
\vskip 5pt

\subheading{\bf References}

\vskip 5pt

[1] Artin M, Grothendieck A, Verdier J.,  SGA 4 , (1963-64),
Lecture Notes in Mathematics 269 and 270, Springer, (1972).

\vskip 3pt

[2] Barr, M., Grillet P.A., ``Exact Categories'' and ``Regular 
Categories'', Springer Lecture Notes 236 (1971).

\vskip 3pt

[3] Dubuc, E.J., ``Adjoint Triangles'', in Springer Lecture Notes 61 
(1968).

\vskip 3pt

[4] Edwards, H. M., ``Galois Theory'', Graduate Texts in Mathematics 
101, Springer (1984).

\vskip 3pt

[5] Freyd, P., ``Abelian Categories'', Harper and Row (1964).

\vskip 3pt

[6] Grothendieck, A., SGA1 (1960-61), Springer Lecture Notes 224 (1971).

\vskip 3pt

[7] Grothendieck, A., ``La Longe Marche a travers la theorie de 
Galois'' (1980-81), edited by Jean Malgoire, Universite Montpellier 
II (1995).

\vskip 3pt

[8] Johnstone, P. T., ``Topos Theory'', Academic Press, (1977).

\vskip 3pt

[9] Joyal, A. Tierney, M.  ``An extension of the Galois Theory of 
Grothendieck''. Memoirs of the American Mathematical Society 151 (1984).

\vskip 3pt

[10] Moerdijk, I. ``Continuous Fibrations and Inverse Limits of 
Toposes'', Compositio Math. 58 (1986).   

\vskip 3pt

[11] Munkres, J. R., ``Topology'', Prentice Hall (1975).

\end